\documentclass[utf8]{frontiersinFPHY_FAMS_nologo} % Vancouver Reference Style (Numbered) for articles in the journals "Frontiers in Physics" and "Frontiers in Applied Mathematics and Statistics"

\setcitestyle{square} % for articles in the journals "Frontiers in Physics" and "Frontiers in Applied Mathematics and Statistics"

\usepackage{url,hyperref,microtype}
\usepackage[onehalfspacing]{setspace}

\usepackage[utf8]{inputenc}
\usepackage{verbatimbox}
\usepackage{diagbox} % for tables

\usepackage{amsmath, amsfonts, amssymb, mathdots}
\usepackage{amsthm} % remove for SISC
\usepackage{hyperref}
\hypersetup{colorlinks,citecolor=blue}

\usepackage{hypcap}
\usepackage{stmaryrd} % for jump

\usepackage{graphicx, color, bm}
\usepackage{caption,subcaption}
\usepackage{booktabs}

\usepackage{tikz, pgfplots, pgfplotstable}
\pgfplotsset{compat=1.16}
\definecolor{markercolor}{RGB}{124.9, 255, 160.65}

\newtheorem{remark}{Remark}

\renewcommand{\hat}{\widehat}
\renewcommand{\tilde}{\widetilde}

\newcommand*\diff[1]{\mathop{}\!{\mathrm{d}#1}} % d in integrand
\newcommand{\diag}[1]{{\rm diag}\LRp{#1}}
\newcommand{\td}[2]{\frac{{\rm d}#1}{{\rm d}{ {#2}}}}
\newcommand{\pd}[2]{\frac{\partial#1}{\partial#2}}
\newcommand{\nor}[1]{\left\| #1 \right\|}
\newcommand{\LRp}[1]{\left( #1 \right)}

\newcommand{\LRb}[1]{\left| #1 \right|}

\newcommand{\fnt}[1]{\bm{\mathsf{ #1}}}

\newcommand{\eqlab}[1]{\begin{align}#1\end{align}}
\newcommand{\bmat}[1]{\begin{bmatrix}#1\end{bmatrix}}

\newcommand{\blue}[1]{\textcolor{blue}{#1}}
\newcommand{\red}[1]{\textcolor{red}{#1}}

\newcommand{\mysubtable}[2]{
{#2}
\\
\vspace{.5em}
\text{\small #1}
}

\graphicspath{{figs/}}

\usepackage{lineno}
\newcommand*\patchAmsMathEnvironmentForLineno[1]{%
  \expandafter\let\csname old#1\expandafter\endcsname\csname #1\endcsname
  \expandafter\let\csname oldend#1\expandafter\endcsname\csname end#1\endcsname
  \renewenvironment{#1}%
     {\linenomath\csname old#1\endcsname}%
     {\csname oldend#1\endcsname\endlinenomath}}%
\newcommand*\patchBothAmsMathEnvironmentsForLineno[1]{%
  \patchAmsMathEnvironmentForLineno{#1}%
  \patchAmsMathEnvironmentForLineno{#1*}}%
\AtBeginDocument{%
\patchBothAmsMathEnvironmentsForLineno{equation}%
\patchBothAmsMathEnvironmentsForLineno{align}%
\patchBothAmsMathEnvironmentsForLineno{flalign}%
\patchBothAmsMathEnvironmentsForLineno{alignat}%
\patchBothAmsMathEnvironmentsForLineno{gather}%
\patchBothAmsMathEnvironmentsForLineno{multline}%
}
\def\keyFont{\fontsize{8}{11}\helveticabold }
\def\firstAuthorLast{Chan {et~al.}} %use et al only if is more than 1 author
\def\Authors{Jesse Chan\,$^{1,*}$, Hendrik Ranocha\,$^{2}$, Andr\'{e}s M. Rueda-Ram\'{i}rez\,$^{3,4}$, Gregor Gassner\,$^{3,4}$, and Tim Warburton\,$^{5}$}
\def\Address{
$^{1}$Department of Computational and Applied Mathematics, Rice University, Houston, Texas, USA\\
$^{2}$Department of Mathematics, University of Hamburg, Hamburg, Germany\\
$^{3}$Department of Mathematics and Computer Science, University of Cologne, Cologne, Germany\\
$^{4}$Center for Data and Simulation Science, University of Cologne, Cologne, Germany\\
$^{5}$Department of Mathematics, Virginia Tech, Blacksburg, Virginia, USA
}
\def\corrAuthor{Jesse Chan}
\def\corrEmail{jesse.chan@rice.edu}

\begin{document}
\onecolumn
\firstpage{1}

\title[On the entropy projection and robustness of entropy stable DG]{On the entropy projection and the robustness of high order entropy stable discontinuous Galerkin schemes for under-resolved flows}

\author[\firstAuthorLast ]{\Authors} %This field will be automatically populated
\address{} %This field will be automatically populated
\correspondance{} %This field will be automatically populated

\extraAuth{}% If there are more than 1 corresponding author, comment this line and uncomment the next one.
%\extraAuth{corresponding Author2 \\ Laboratory X2, Institute X2, Department X2, Organization X2, Street X2, City X2 , State XX2 (only USA, Canada and Australia), Zip Code2, X2 Country X2, email2@uni2.edu}

\maketitle

\begin{abstract}
High order entropy stable schemes provide improved robustness for computational simulations of fluid flows. However, additional stabilization and positivity preserving limiting can still be required for variable-density flows with under-resolved features. We demonstrate numerically that entropy stable DG methods which incorporate an ``entropy projection'' are less likely to require additional limiting to retain positivity for certain types of flows. We conclude by investigating potential explanations for this observed improvement in robustness.

\tiny
 \keyFont{ \section{Keywords:} computational fluid dynamics, high order, discontinuous Galerkin, summation-by-parts, entropy stability, robustness} %All article types: you may provide up to 8 keywords; at least 5 are mandatory.
\end{abstract}

\section{Introduction}

Discontinuous Galerkin (DG) schemes have received interest within computational fluid dynamics (CFD) due to their high order accuracy and ability to handle unstructured curved meshes. In particular, there has been interest in DG methods for simulations of under-resolved flows \cite{gassner2013accuracy, beck2014high, moura2017eddy, fernandez2018ability, lv2018underresolved}. Among such schemes, ``entropy stable'' DG methods based on a ``flux differencing'' formulation have received interest due to their robustness with respect to shocks and turbulence \cite{flad2017use, winters2018comparative, rojas2021robustness, parsani2021high}.

Entropy conservative and entropy stable flux differencing schemes were originally formulated for finite difference methods in \cite{fjordholm2012arbitrarily, fisher2013high}. They were extended to tensor product grids using discontinuous spectral collocation schemes (also known as discontinuous Galerkin spectral element methods, or DGSEM) \cite{carpenter2014entropy, gassner2016split}. Entropy stable collocation schemes were extended to simplicial meshes in \cite{chen2017entropy, crean2018entropy} using multi-dimensional summation-by-parts (SBP) operators \cite{hicken2016multidimensional}. Non-collocation entropy stable schemes have also been developed. These schemes began with staggered grid schemes on tensor product grids in \cite{parsani2016entropy}, which were later extended to simplicial elements in \cite{fernandez2019staggered}. ``Modal'' entropy stable DG formulations \cite{chan2017discretely, chan2018discretely, chan2019skew} have been utilized to construct a variety of new entropy stable schemes, including Gauss DG methods \cite{chan2018efficient, chan2020mortar} and reduced order models \cite{chan2019entropy}. We note that under appropriate choices of quadrature, these ``modal'' formulations reduce to collocation-type entropy stable schemes. Entropy stable schemes have since been extended to an even wider array of discretizations, such as line DG methods, discontinuous Galerkin difference methods, and $C^0$ continuous discretizations \cite{pazner2019analysis, hicken2020entropy, yan2021entropy}.

The main difference between non-collocation and collocation-type entropy stable schemes is the use of transformations between conservative variables and entropy variables together with projection or prolongation operators to facilitate a discrete proof of entropy stability. This is referred to as the ``entropy projection'' in \cite{chan2017discretely, pazner2019analysis} and as the interpolation or prolongation of entropy variables in \cite{parsani2016entropy, yan2021entropy}. This approach is also equivalent to the mixed formulation of \cite{gkanis2021new}. We will refer to this transformation as the ``entropy projection'' for the remainder of the paper.

The motivation for introducing the entropy projection has been to enable the use of more accurate quadrature rules or novel basis functions. This has been at the cost of additional complexity and issues related to the sensitivity of the entropy variables for near-vacuum states \cite{chan2017discretely, yan2021entropy}. To the best of the authors' knowledge, no inherent advantages in using the entropy projection have been observed in the literature. This paper focuses on the following observation: high order entropy stable schemes based on the entropy projection appear to be more robust than entropy stable collocation schemes for two and three dimensional simulations of under-resolved variable-density fluid flows with small-scale features.

The structure of the paper is as follows: Section~\ref{sec:formulation} reviews mathematical formulations of entropy stable schemes which involve the entropy projection. Section~\ref{sec:comparison} documents the observed difference in robustness for a variety of problems in two and three dimensions, and provides analysis and numerical experiments which support that the primary difference between unstable and stable schemes is the entropy projection. Section~\ref{sec:conjecture} conjectures potential explanations for why the entropy projection might improve robustness. We conclude with Section~\ref{sec:spectra}, which explores potential applications towards under-resolved flow simulations.

\section{Formulation of high order entropy stable DG schemes}
\label{sec:formulation}

In this section, we provide a brief description of high order entropy stable schemes in 1D. More detailed derivations, multi-dimensional formulations, and extensions to curved grids can be found in \cite{chen2017entropy, chan2017discretely, chan2018efficient, crean2018entropy, chan2019entropy, chan2019skew}.

The notation in this paper is motivated by notation in \cite{crean2018entropy, fernandez2019entropy}. Unless otherwise specified, vector and matrix quantities are denoted using lower and upper case bold font, respectively.  Spatially discrete quantities are denoted using a bold sans serif font. Finally, the output of continuous functions evaluated over discrete vectors is interpreted as a discrete vector.

For example, if $\fnt{x}$ denotes a vector of point locations, i.e., $(\fnt{x})_i = \bm{x}_i$, then $u(\fnt{x})$ is interpreted as the vector
\[	({u}(\fnt{x}))_i = {u}(\bm{x}_i).
\]
Similarly, if $\fnt{u} = {u}(\fnt{x})$, then ${f}(\fnt{u})$ corresponds to the vector
\[
	({f}(\fnt{u}))_i = {f}(u(\bm{x}_i)).
\]
Vector-valued functions are treated similarly. For example, given a vector-valued function $\bm{f}:\mathbb{R}^n\rightarrow \mathbb{R}^n$ and a vector of coordinates $\fnt{x}$, we adopt the convention that $\LRp{\bm{f}(\fnt{x})}_i = \bm{f}(\bm{x}_i)$.

\subsection{Conservation laws with entropy}

In this section, we review the construction of entropy conservative and entropy stable schemes for a one-dimensional system of nonlinear conservation laws
\[
\pd{\bm{u}}{t} + \pd{\bm{f}(\bm{u})}{x} = \bm{s}(\bm{u}),
\]
where $\bm{s}(\bm{u})$ is a source term.  We assume the domain is exactly represented by a uniform mesh consisting of non-overlapping intervals $D^k$, and that the solution $u(x)$ is approximated by degree $N$ polynomials over each element. We also introduce entropy conservative numerical fluxes $\bm{f}_{S}(\bm{u}_L,\bm{u}_R)$ \cite{tadmor1987numerical}, which are bivariate functions of ``left'' and ``right'' states $\bm{u}_L, \bm{u}_R$. In addition to being symmetric and consistent, entropy conservative numerical fluxes satisfy an ``entropy conservation'' property
\eqlab{
\LRp{\bm{v}_L-\bm{v}_R}^T\bm{f}_{S}(\bm{u}_L,\bm{u}_R) = \psi(\bm{u}_L) - \psi(\bm{u}_R).
\label{eq:ecflux}
}
Here, $\bm{v}_L, \bm{v}_R$ are entropy variables evaluated at the left and right states, and $\psi(\bm{u})$ denotes the ``entropy potential''. Examples of expressions for entropy variables and entropy potentials can be found in \cite{chen2017entropy}.

\subsection{Collocation formulations}

Degree $N$ entropy stable collocation schemes are typically built from Legendre-Gauss-Lobatto (LGL) quadrature rules with $(N+1)$ points. Let $\fnt{x}, \fnt{w}$ denote vectors of quadrature points and weights on the reference interval $[-1,1]$. Let $\ell_i(x)$ denote Lagrange polynomials at LGL nodes, and let $\fnt{u}$ denote the vector of solution nodal values $u(x_i)$. Define the matrices
\begin{gather*}
\fnt{M} = \diag{\fnt{w}}, \qquad \fnt{Q}_{ij} = \int_{-1}^1 \pd{\ell_j}{x}\ell_i \diff{x},
\\
\fnt{B} = \begin{bmatrix}
-1 &\\
& 1
\end{bmatrix}, \qquad \fnt{V}_f = \begin{bmatrix}
1 & \ldots & 0\\
0 & \ldots & 1
\end{bmatrix}.
\end{gather*}
Here $\fnt{V}_f$ is a face interpolation or extraction matrix which maps from volume nodes to face nodes.
Flux derivatives are discretized using a ``flux differencing'' approach involving summation-by-parts (SBP) operators and entropy conservative fluxes \cite{tadmor1987numerical}. An entropy stable collocation formulation can now be defined on an element $D^k$ as follows:
\begin{align} \label{eq:collocationDisc}
h \fnt{M} \td{\fnt{u}}{t} + \LRp{\LRp{\fnt{Q}-\fnt{Q}^T}\circ\fnt{F}}\fnt{1} + \fnt{V}_f^T\fnt{B}\fnt{f}^* = \fnt{s}(\fnt{u}), \qquad \fnt{F}_{ij} = \bm{f}_S\LRp{\fnt{u}_i, \fnt{u}_j},
\end{align}
where $h$ is the size of the element $D^k$ and $\circ$ denotes the matrix Hadamard product {\cite{fjordholm2012arbitrarily, fisher2013high, carpenter2014entropy}}.\footnote{Since the entries of $\fnt{F}$ are vector-valued, the Hadamard product $\LRp{\fnt{Q}-\fnt{Q}^T}\circ\fnt{F}$ should be understood as each scalar entry of $\LRp{\fnt{Q}-\fnt{Q}^T}$ multiplying each vector-valued entry of ${\fnt{F}}$.}
Here, $\fnt{f}^*$ is a vector which contains numerical fluxes at the left and right endpoints of the interval
\[
\fnt{f}^* = \begin{bmatrix}
\bm{f}^*(\fnt{u}_1^+, \fnt{u}_1)\\
\bm{f}^*(\fnt{u}_{N+1}, \fnt{u}_{N+1}^+)
\end{bmatrix},
\]
where $\fnt{u}_1^+, \fnt{u}_{N+1}^+$ denote exterior nodal values on neighboring elements. If $\fnt{f}^*$ is an entropy conservative flux, then the resulting numerical method is semi-discretely entropy conservative. If $\fnt{f}^*$ is an entropy stable flux (for example, Lax-Friedrichs flux, HLLC, and certain matrix penalizations \cite{chen2017entropy, winters2017uniquely}) then the resulting scheme also dissipates entropy.

\subsection{``Modal'' formulations}

Degree $N$ entropy stable ``modal'' DG schemes generalize collocation schemes to arbitrary choices of quadrature. In one dimension, this allow for the use of higher accuracy volume quadratures. In higher dimensions, modal formulations also enable more general choices of surface quadrature. These schemes introduce an additional ``entropy projection'' step to facilitate the semi-discrete proof of entropy stability or conservation.

We now assume the solution is represented using some arbitrary basis over each element, such that $u(x) = \sum_j \fnt{u}_j \phi_j(x)$. Let $\bm{x}, \bm{w}$ now denote a general quadrature rule with positive quadrature weights. We define quadrature-based interpolation matrices $\fnt{V}_q, \fnt{V}_f$, the mass matrix $\fnt{M}$, and the modal differentiation matrix $\hat{\fnt{Q}}$
\begin{gather*}
\LRp{\fnt{V}_q}_{ij} = \phi_j(\fnt{x}_i), \qquad \LRp{\fnt{V}_f}_{1j} = \phi_j\LRp{-1}, \qquad \LRp{\fnt{V}_f}_{2j} = \phi_j\LRp{1}, \\
\fnt{M} = \fnt{V}_q^T\diag{\fnt{w}}\fnt{V}_q,  \qquad \hat{\fnt{Q}}_{ij} = \int_{-1}^1 \pd{\phi_j}{x} \phi_i \diff{x}.
\end{gather*}
We introduce the quadrature-based projection matrix $\fnt{P}_q= \fnt{M}^{-1}\fnt{V}_q^T\diag{\fnt{w}}$. Using $\fnt{P}_q$ and $\hat{\fnt{Q}}$, we can construct quadrature-based differentiation and extrapolation matrices {$\fnt{Q}$}, $\fnt{E}$
\[
\fnt{Q} = \fnt{P}_q^T\hat{\fnt{Q}}\fnt{P}_q, \qquad \fnt{E} = \fnt{V}_f\fnt{P}_q.
\]
To accommodate general quadrature rules which may not include boundary points, we introduce hybridized SBP operators $\fnt{Q}_{h}$ on the reference interval $[-1,1]$
\[
\fnt{Q}_{h} = \frac{1}{2}\bmat{
\fnt{Q} - \LRp{\fnt{Q}}^T & \fnt{E}^T\fnt{B}\\
-\fnt{B}\fnt{E}  & \fnt{B}
}.
\]
The use of such operators simplifies the implementation for general quadrature rules and nodal sets which do not include boundary nodes \cite{chan2017discretely, chen2020review}.
Next, we define $\fnt{V}_h$ as the interpolation matrix to \textit{both} volume and surface quadrature points
\[
\fnt{V}_h = \bmat{
\fnt{V}_q\\
\fnt{V}_f
}.
\]
We also introduce the $L^2$ projection of the entropy variables and the ``entropy projected'' conservative variables $\tilde{\fnt{u}}$
\[
\fnt{v} = \fnt{P}_q \bm{v}\LRp{\fnt{V}_q\fnt{u}}, \qquad \tilde{\fnt{u}} = \bm{u}\LRp{\fnt{V}_h\fnt{v}},
\]
which are defined by evaluating the mapping from entropy to conservative variables using the projected entropy variables. Note that the projected entropy variables $\fnt{v}$ is a vector corresponding to modal coefficients, while $\tilde{\fnt{u}}$ corresponds to point values at volume and face quadrature points.

An entropy stable modal DG discretization over a single element $D^k$ is then
\begin{gather} \label{eq:modalDisc}
h\fnt{M}\td{\fnt{u}}{t} + \fnt{V}_h^T\LRp{\LRp{\fnt{Q}_{h} - \fnt{Q}_{h}^T} \circ \fnt{F}}\fnt{1} + \fnt{V}_f^T \fnt{B} \fnt{f}^*= \fnt{s}(\fnt{u}),
\\
\fnt{F}_{ij} = \bm{f}_S\LRp{\tilde{\fnt{u}}_i,\tilde{\fnt{u}}_j}, \qquad \fnt{f}^* = \begin{bmatrix}
\fnt{f}^*\LRp{\tilde{\fnt{u}}_1^+,\tilde{\fnt{u}}_1}  \\
\fnt{f}^*\LRp{\tilde{\fnt{u}}_{N+1},\tilde{\fnt{u}}_{N+1}^+}
\end{bmatrix}. \nonumber
\end{gather}
Note that the right hand side formulation is evaluated not using the conservative variables $\fnt{u}$, but the ``entropy projected'' conservative variables $\tilde{\fnt{u}}$.

While we have presented entropy stable DG schemes using a general ``modal'' DG framework, the formulation reduces to existing methods under appropriate choices of quadrature and basis. For example, specifying LGL quadrature on a tensor product element recovers entropy stable spectral collocation schemes \cite{chan2018efficient}. SBP discretizations without an underlying basis on simplices \cite{hicken2016multidimensional, chen2017entropy, crean2018entropy} can also be recovered for appropriate quadrature rules by redefining the interpolation and projection matrices $\fnt{V}_q, \fnt{P}_q$ \cite{wu2020high}.

\section{Numerical comparisons of collocation and entropy projection schemes}
\label{sec:comparison}

In this section, we will demonstrate numerically that a significant difference in robustness is observed between collocation and entropy projection-based discretizations of the Euler and ideal MHD equations.
For the Euler equations, we study the Kelvin-Helmholtz, Rayleigh-Taylor, and Richtmeyer-Meshkov instabilities, and for the MHD equations we study a magnetized Kelvin-Helmholtz instability. 
All of these examples exhibit small-scale turbulent-like features. Moreover, we observe a difference in robustness between entropy stable collocation and entropy projection-based methods independently of the polynomial degrees, mesh resolutions, and type of mesh (e.g., quadrilateral or triangular). We focus on the following entropy stable DG methods:
\begin{itemize}
\item On quadrilateral meshes:
\begin{enumerate}
\item DGSEM: collocation scheme based on the tensor product of one-dimensional $(N+1)$ point LGL quadrature,
\item Gauss DG: a ``collocation'' scheme based on the tensor product of one-dimensional $(N+1)$ point Gauss quadrature. The entropy projection is used to evaluate interface fluxes \cite{chan2018efficient},
\end{enumerate}
\item On triangular meshes:
\begin{enumerate}
\item SBP: a collocation scheme based on multi-dimensional summation-by-parts finite difference operators \cite{hicken2016multidimensional, chen2017entropy},
\item Modal: a modal formulation utilizing quadrature rules which exactly integrate entries of the volume and face mass matrices \cite{chan2017discretely}.
\end{enumerate}
\end{itemize}

\begin{remark}
It is known that the Kelvin-Helmholtz, Rayleigh-Taylor, and Richtmeyer-Meshkov instabilities are notoriously sensitive to initial conditions and discretization parameters, and that numerical schemes may not converge to a unique solution \cite{fjordholm2017construction, schroeder2019reference}. Instead, this paper focuses on these problems as stress tests of robustness.
\end{remark}

Unless specified otherwise, all numerical experiments utilize a Lax-Friedrichs interface flux with Davis wavespeed estimate \cite{davis1988simplified}. 
We also experimented with HLL and HLLC surface fluxes, but did not notice a significant difference. We also note that instead of discontinuous initial conditions, we utilize smoothed approximations for each problem considered here.

All experiments are also performed on uniform meshes. For triangular meshes, this mesh is constructed by bisecting each element of a uniform quadrilateral mesh along the diagonal. Unless specified otherwise, all results are produced using the Julia \cite{bezanson2017julia} simulation framework Trixi.jl \cite{schlottkelakemper2021purely, ranocha2022adaptive}. For most experiments, we utilize an optimized adaptive 4th order 9-stage Runge-Kutta method \cite{ranocha2021optimized} implemented in OrdinaryDiffEq.jl \cite{rackauckas2017differentialequations}. The absolute and relative tolerances are set to $10^{-7}$ unless specified otherwise. Scripts generating main results are included in a companion repository for reproducibility \cite{reprorepo}.

We note that the robustness, efficiency, and high order accuracy of both entropy stable DGSEM and entropy stable Gauss DG schemes have been verified in previous works \cite{winters2018comparative, chan2018efficient, rojas2021robustness, chan2020mortar, parsani2021high}, and will not be addressed in detail in this paper. However, the difference in robustness between the two methods has not been previously observed in the literature, and will be the focus of this work.

\subsection{Euler Equations of Gas Dynamics}

We consider first the two and three-dimensional problems for the Euler equations of gas dynamics.
The conservative variables for the three-dimensional Euler equations are density, momentum, and total energy, $\bm{u} = (\rho, \rho \fnt{v}, E)$, where the vector $\fnt{v} = (u,v,w)$ contains the velocities in $x$, $y$ and $z$, respectively.
The flux reads
\[
\bm{f}(\bm{u}) = 
\begin{pmatrix}
\rho \fnt{v} \\
\rho \fnt{v} \fnt{v}^T + \fnt{I} p \\
\fnt{v} \left( \frac{1}{2} \rho \|\fnt{v}\|^2 + \frac{\gamma p}{\gamma-1} \right)
\end{pmatrix},
\]
where $\fnt{I}$ is the $3 \times 3$ identity matrix, $\gamma$ is the heat capacity ratio, and $p = (\gamma-1)\left(E - \rho \|\fnt{v}\|^2/2 \right)$ is the gas pressure. For two-dimensional problems, we neglect the third component of the velocity, $w$, and $\fnt{I}$ becomes the $2 \times 2$ identity matrix.

All the following experiments use the entropy conservative and kinetic energy preserving flux of Ranocha \cite{ranocha2020entropy, ranocha2021preventing}; however, similar results were observed when experimenting with the entropy conservative flux of Chandrashekar \cite{chandrashekar2013kinetic}.

%\subsection{Two-dimensional experiments}

\subsubsection{Two dimensional Kelvin-Helmholtz instability} \label{eq:KHI_euler}

\begin{figure}
\centering
\includegraphics[width=\textwidth]{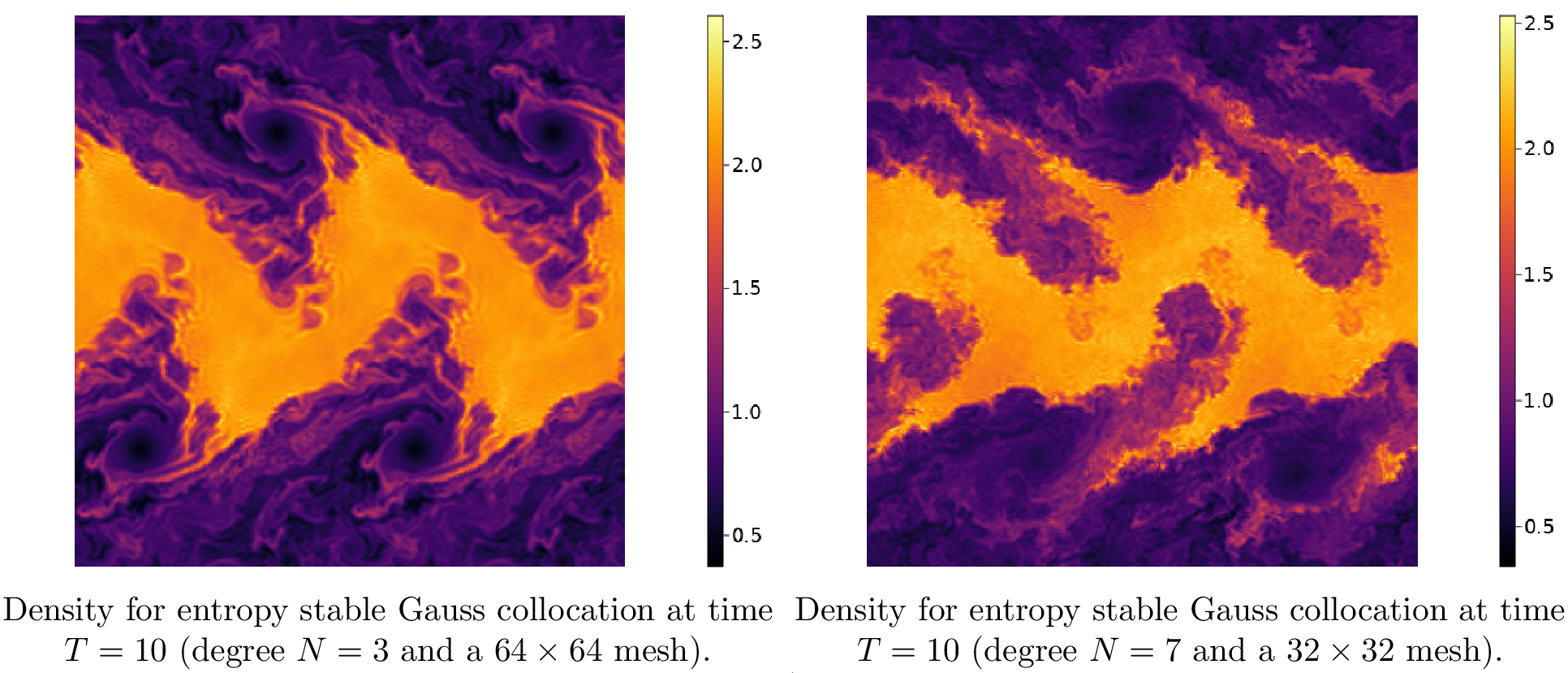}
%\subcaptionbox*{Density for entropy stable Gauss DG at time $T = 10$ (degree $N=3$ and a $64\times 64$ mesh).}{\includegraphics[width=.49\textwidth]{khi_p3_64_cells.png}}
%\hspace{.25em}
%\subcaptionbox*{Density for entropy stable Gauss DG at time $T = 10$ (degree $N=7$ and a $32\times 32$ mesh).}{\includegraphics[width=.49\textwidth]{khi_p7_32_cells.png}}
\caption{Snapshots of density for the Kelvin-Helmholtz instability using an entropy stable Gauss DG scheme on uniform quadrilateral meshes.}
\label{fig:khi_intro}
\end{figure}

We perform additional experiments analyzing the robustness of entropy stable DGSEM and Gauss DG for the Kelvin-Helmholtz instability. The domain is $[-1,1]^2$ with initial condition from \cite{rueda2021subcell}:
\begin{align}
&\rho = \frac{1}{2} + \frac{3}{4}B,&
%&\rho = \begin{cases}
%2 & y \in [-1/2, 1/2]\\
%1/2 & \text{otherwise}
%\end{cases}
&p = 1,& \nonumber\\
&u = \frac{1}{2}(B-1),&
&v = \frac{1}{10}\sin\LRp{2\pi x},& \label{eq:khi_ic}
\end{align}
where $B(x,y)$ is a smoothed approximation to a discontinuous step function
\begin{equation}
B(x,y) = \tanh(15 y + 7.5) - \tanh(15 y - 7.5).
\label{eq:B}
\end{equation}
Each solver is run until final time $T_{\rm final} = 15$. 
As can be observed in Figure~\ref{fig:khi_intro}, the solution differs significantly between the $N=3$ and $N=7$ simulations. This is likely a consequence of the well-known sensitivity of the Kelvin-Helmholtz instability to small perturbations and numerical resolutions \cite{fjordholm2017construction, schroeder2019reference}.

\begin{table}
\centering
\mysubtable{KHI, quadrilateral mesh, $N_{\rm cells} = 16$}{
\centering
\begin{tabular}{|c||c|c|c|c|c|c|c|}
\hline
\diagbox{Solver}{Degree} & 1 & 2 & 3 & 4 & 5 & 6 & 7\\
\hline
{Collocation} &\blue{15}  & \red{4.807}  & \red{3.769}  & \red{4.433}  & \red{3.737}  & \red{3.369}  & \red{3.642}  \\
\hline
{Entropy projection} & \blue{15} & \blue{15} & \blue{15} & \blue{15} & \blue{15} & \blue{15} & \blue{15}  \\
\hline
\end{tabular}
}\\
\vspace{1em}
\mysubtable{KHI, quadrilateral mesh, $N_{\rm cells} = 32$}{
\centering
\begin{tabular}{|c||c|c|c|c|c|c|c|}
\hline
\diagbox{Solver}{Degree} & 1 & 2 & 3 & 4 & 5 & 6 & 7\\
\hline
{Collocation} & \blue{15}  & \red{4.116}  & \red{3.652}  & \red{4.266}  & \red{3.54}  & \red{3.663}  & \red{3.556}  \\
\hline
{Entropy projection} & \blue{15} & \blue{15} & \blue{15} & \blue{15} & \blue{15} & \blue{15} & \blue{15}  \\
\hline
\end{tabular}
}\\
\vspace{1em}
\mysubtable{KHI, triangular mesh, $N_{\rm cells} = 16$}{
\centering
\begin{tabular}{|c||c|c|c|c|c|c|}
\hline
\diagbox{Solver}{Degree} & 1 & 2 & 3 & 4 & 5 & 6 \\
\hline
{Collocation} & \blue{15}  & \red{3.984}  & \red{3.441}  & \red{2.993}  & \red{2.943}  & \red{3.128}  \\
\hline
{Entropy projection} & \blue{15} & \blue{15} & \blue{15} & \blue{15} & \blue{15} & \blue{15}  \\
\hline
\end{tabular}
}\\
\vspace{1em}
\mysubtable{KHI, triangular mesh, $N_{\rm cells} = 32$}{
\centering
\begin{tabular}{|c||c|c|c|c|c|c|}
\hline
\diagbox{Solver}{Degree} & 1 & 2 & 3 & 4 & 5 & 6 \\
\hline
{Collocation} & \red{3.919}  & \red{3.452}  & \red{3.191}  & \red{2.958}  & \red{3.063}  & \red{3.269}  \\
\hline
{Entropy projection} & \blue{15}  & \blue{15}  & \blue{15}  & \blue{15}  & \blue{15}  & \blue{15}  \\
\hline
\end{tabular}
}\\
\vspace{.5em}
\caption{End time for simulations of the Kelvin-Helmholtz instability on quadrilateral and triangular meshes. On quadrilateral meshes, ``collocation'' refers to a nodal DGSEM discretization, while ``entropy projection'' refers to a method based on Gauss nodes.  On triangular meshes, ``collocation'' refers to nodal SBP discretization, while ``entropy projection'' refers to a modal entropy stable DG method.  Times colored \blue{blue} correspond to simulations which did not crash and ran to completion, while times colored \textcolor{red}{red} denote simulations which did crash.}
\label{tab:compare_khi}
\end{table}

\subsubsection{Two dimensional Rayleigh-Taylor instability}  The two-dimensional Rayleigh-Taylor instability generates small-scale flow features through buoyancy or gravity effects \cite{richtmyer1954taylor, youngs1984numerical}. The setup involves a heavy and light fluid suspended above one another separated by a curved interface, and buoyancy or gravity results in displacement of the lighter fluid into the heavier one. This displacement causes velocity shear and the formation of additional Kelvin-Helmholtz instabilities along the interface. The domain is $[0, 1/4] \times [0, 1]$.

\begin{figure}
\centering
\includegraphics[width=\textwidth]{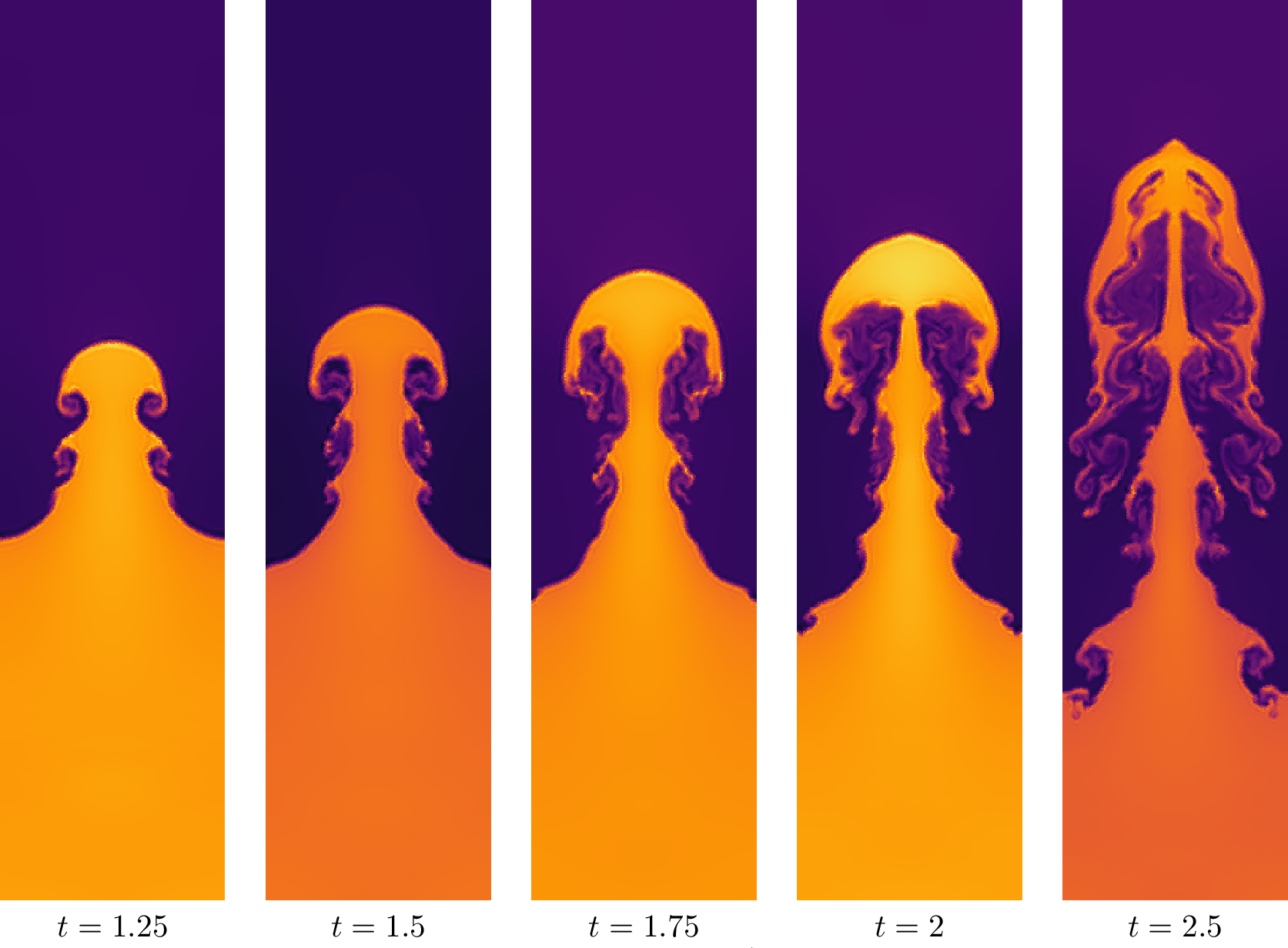}
%\subcaptionbox*{$t = 1.25$}{\includegraphics[width=.15\textwidth]{rti_density_Gauss_p3_32_cells_t1p25.png}}
%\hspace{.6em}
%\subcaptionbox*{$t = 1.5$}{\includegraphics[width=.15\textwidth]{rti_density_Gauss_p3_32_cells_t1p5.png}}
%\hspace{.6em}
%\subcaptionbox*{$t = 1.75$}{\includegraphics[width=.15\textwidth]{rti_density_Gauss_p3_32_cells_t1p75.png}}
%\hspace{.6em}
%\subcaptionbox*{$t = 2$}{\includegraphics[width=.15\textwidth]{rti_density_Gauss_p3_32_cells_t2.png}}
%\hspace{.6em}
%\subcaptionbox*{$t = 2.5$}{\includegraphics[width=.15\textwidth]{rti_density_Gauss_p3_32_cells_t2p5.png}}
\caption{Density for a Rayleigh-Taylor instability for a degree $N=3$ entropy stable Gauss DG scheme on a mesh of $32 \times 128$ elements.}
\label{fig:rti}
\end{figure}

Let $d_{a,b}(x) = a + \frac{1}{2} (1 + \tanh(s x)) (b - a)$ denote a smoothed approximation (with slope $s$) to a discontinuous function with values $a$ for $x < 0$ and $b$ for $x > 0$.
The initial condition is given by
\begin{align*}
%&\rho = \begin{cases}
%1 & y < 1/2\\
%2 & y > 1/2
%\end{cases},&
&\rho = d_{2, 1}\LRp{y - \frac{1}{2}},&
&p = \begin{cases}
2y +1 & y < 1/2\\
y + 3/2 & y \geq 1/2
\end{cases},&\\
&u = 0&
&v = -\frac{c}{40} \cos(8 k\pi x) \sin(\pi y)^6,&
\end{align*}
where $c = \sqrt{\gamma p / \rho}$ is the speed of sound. Here, we borrow from \cite{remacle2003adaptive} and multiply the $y$-velocity perturbation by $\sin(\pi y)^6$ so that $u, v$ satisfy wall boundary conditions. We also add gravity source terms to the $y$-momentum and energy equations:
\begin{align*}
\bm{s}(\bm{x}, t) = \begin{bmatrix}
0 & 0 & g \rho & g \rho v
\end{bmatrix},
\end{align*}
where $s=15$. Note that the sign of gravity is such that the light fluid flows up into the heavy fluid.
Reflective wall boundary conditions are imposed at all boundaries using mirror states, which results in an entropy stable scheme under the Lax-Friedrichs flux \cite{chen2017entropy, hindenlang2019stability}. Figure~\ref{fig:rti} shows snapshots of the density for a degree $N=3$ entropy stable Gauss DG scheme on a mesh of $32 \times 128$ elements at various times.

\begin{table}
\centering
\mysubtable{RTI, quadrilateral mesh, $N_{\rm cells} = 16$}{
\centering
\begin{tabular}{|c||c|c|c|c|c|c|c|}
\hline
\diagbox{Solver}{Degree} & 1 & 2 & 3 & 4 & 5 & 6 & 7\\
\hline
{Collocation} & \red{3.674}  & \red{3.44}  & \red{3.332}  & \red{3.257}  & \red{3.106}  & \red{3.034}  & \red{3.044}  \\
\hline
{Entropy projection} & \blue{15} & \blue{15} & \blue{15} & \blue{15} & \blue{15} & \blue{15} & \blue{15}  \\
\hline
\end{tabular}
}\\
\vspace{1em}
\mysubtable{RTI, quadrilateral mesh, $N_{\rm cells} = 32$}{
\centering
\begin{tabular}{|c||c|c|c|c|c|c|c|}
\hline
\diagbox{Solver}{Degree} & 1 & 2 & 3 & 4 & 5 & 6 & 7\\
\hline
{Collocation} & \red{3.996}  & \red{3.144}  & \red{3.44}  & \red{3.155}  & \red{3.031}  & \red{2.972}  & \red{2.976}  \\
\hline
{Entropy projection} & \blue{15} & \blue{15} & \blue{15} & \blue{15} & \blue{15} & \blue{15} & \blue{15}  \\
\hline
\end{tabular}
}\\
\vspace{1em}
\mysubtable{RTI, triangular mesh, $N_{\rm cells} = 16$}{
\centering
\begin{tabular}{|c||c|c|c|c|c|c|}
\hline
\diagbox{Solver}{Degree} & 1 & 2 & 3 & 4 & 5 & 6 \\
\hline
{Collocation} & \red{4.297}  & \red{2.87}  & \red{3.238}  & \red{3.229}  & \red{2.927}  & \red{2.881} \\
\hline
{Entropy projection} & \blue{15} & \blue{15} & \blue{15} & \blue{15} & \blue{15} & \blue{15}  \\
\hline
\end{tabular}
}\\
\vspace{1em}
\mysubtable{RTI, triangular mesh, $N_{\rm cells} = 32$}{
\centering
\begin{tabular}{|c||c|c|c|c|c|c|}
\hline
\diagbox{Solver}{Degree} & 1 & 2 & 3 & 4 & 5 & 6 \\
\hline
{Collocation} & \red{3.6}  & \red{2.896}  & \red{3.197}  & \red{3.227}  & \red{3.032}  & \red{2.778}  \\
\hline
{Entropy projection} &\blue{15}  & \blue{15}  & \blue{15}  & \blue{15}  & \blue{15}  & \blue{15}  \\
\hline
\end{tabular}
}\\
\vspace{.5em}
\caption{End time for simulations of the Rayleigh-Taylor instability on quadrilateral and triangular meshes.  On quadrilateral meshes, ``collocation'' refers to a nodal DGSEM discretization, while ``entropy projection'' refers to a method based on Gauss nodes.  On triangular meshes, ``collocation'' refers to nodal SBP discretization, while ``entropy projection'' refers to a modal entropy stable DG method.  Times colored \blue{blue} correspond to simulations which did not crash and ran to completion, while times colored \textcolor{red}{red} denote simulations which did crash.}
\label{tab:compare_rti}
\end{table}

\subsubsection{Two dimensional Richtmeyer-Meshkov instability}

\begin{figure}
\centering
\includegraphics[width=\textwidth]{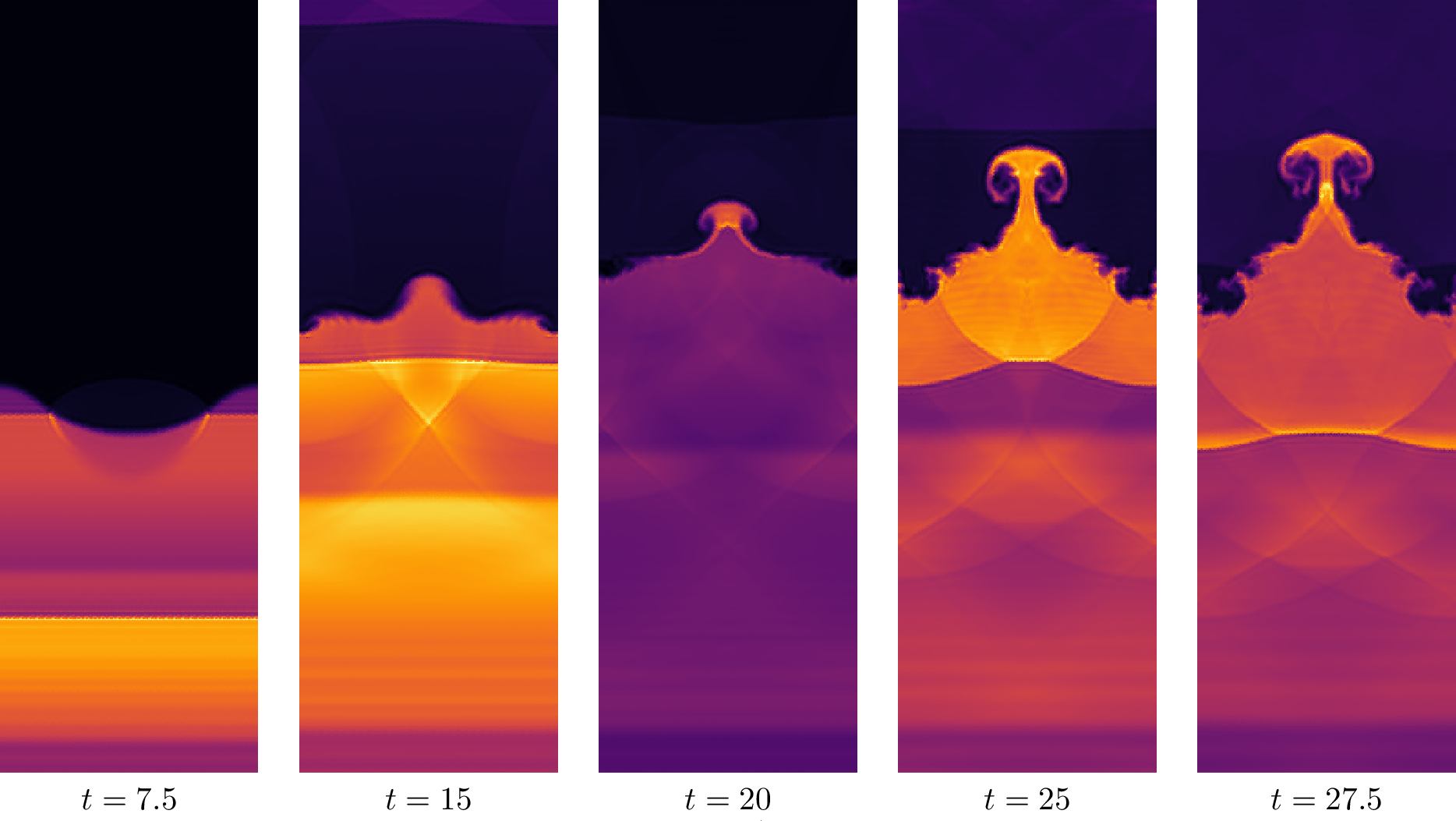}
%\subcaptionbox*{$t = 7.5$}{\includegraphics[width=.17\textwidth]{rmi_density_Gauss_p3_32_cells_t7p5.png}}
%\hspace{.6em}
%\subcaptionbox*{$t = 15$}{\includegraphics[width=.17\textwidth]{rmi_density_Gauss_p3_32_cells_t15.png}}
%\hspace{.6em}
%\subcaptionbox*{$t = 20$}{\includegraphics[width=.17\textwidth]{rmi_density_Gauss_p3_32_cells_t20.png}}
%\hspace{.6em}
%\subcaptionbox*{$t = 25$}{\includegraphics[width=.17\textwidth]{rmi_density_Gauss_p3_32_cells_t25.png}}
%\hspace{.6em}
%\subcaptionbox*{$t = 27.5$}{\includegraphics[width=.17\textwidth]{rmi_density_Gauss_p3_32_cells_t27p5.png}}
\caption{Density for the Richtmeyer-Meshkov instability using a degree $N=3$ entropy stable Gauss DG with $32 \times 96$ elements. The domain is $[0, 40/3] \times [0, 40]$.}
\label{eq:rmi}
\end{figure}

The Richtmeyer-Meshkov instability generates small-scale flow features by passing a shock over a stratified fluid \cite{richtmyer1954taylor, meshkov1969instability}. The domain for this setup is $[0, 40/3] \times [0, 40]$, and the initial density and pressure are given by
\begin{gather*}
\rho = d_{1, \frac{1}{4}}\LRp{y - \LRp{18 + 2 \cos\LRp{\frac{6 \pi x}{L}}}} + d_{3.22, 0}\LRp{\LRb{y-4} - 2},\\
%\rho = smoothed_heaviside(x[2] - (18 + 2 * cos(2 * pi * 3 / L * x[1])), 1.0, .25)
%        \rho = rho + smoothed_heaviside(abs(x[2] - 4) - 2, 3.22, 0.0) # 2 < x < 6
%        p = smoothed_heaviside(abs(x[2] - 4) - 2, 4.9, 1.0)
p = d_{4.9, 1}\LRp{\LRb{y-4} - 2},
%\begin{pmatrix}
%\rho\\
%p
%\end{pmatrix} = \begin{cases}
%\begin{pmatrix}
%1\\
%1
%\end{pmatrix}
%& y \geq 18 + 2 \cos(2 \pi \frac{3x}{40})\\
%\begin{pmatrix}
%4.22\\
%4.9
%\end{pmatrix}
%& 2 \geq y \geq 6 \\
%\begin{pmatrix}
%1/4\\
%1
%\end{pmatrix}
%& \text{otherwise}
%\end{cases}.
\end{gather*}
where we again set the slope $s=15$.
The initial velocities are both set to zero, i.e., $u,v = 0$. We approximate the discontinuous initial condition using smoothed Heaviside functions with a slope of $s=2$ due to the size of the domain.  Reflective wall boundary conditions are imposed everywhere. 
Figure~\ref{eq:rmi} shows pseudocolor plots of the density using a degree $N=3$ entropy stable Gauss DG on a uniform mesh of $32 \times 96$ quadrilateral elements.

\begin{table}
\centering
\mysubtable{RMI, quadrilateral mesh, $N_{\rm cells} = 16$}{
\centering
\begin{tabular}{|c||c|c|c|c|c|c|c|}
\hline
\diagbox{Solver}{Degree} & 1 & 2 & 3 & 4 & 5 & 6 & 7\\
\hline
{Collocation} & \blue{30}  & \blue{30}  & \red{27.96}  & \red{24.94}  & \red{8.851}  & \red{8.853}  & \red{8.85}  \\
\hline
{Entropy projection} &\blue{30}  &\blue{30}  &\blue{30}  &\blue{30}  &\blue{30}  &\blue{30}  &\blue{30}   \\
\hline
\end{tabular}
}\\
\vspace{1em}
\mysubtable{RMI, quadrilateral mesh, $N_{\rm cells} = 32$}{
\centering
\begin{tabular}{|c||c|c|c|c|c|c|c|}
\hline
\diagbox{Solver}{Degree} & 1 & 2 & 3 & 4 & 5 & 6 & 7\\
\hline
{Collocation} & \blue{30}  & \red{25.52}  & \red{23.34}  & \red{8.759}  & \red{7.808}  & \red{7.014}  & \red{7.01}  \\
\hline
{Entropy projection} &\blue{30}  &\blue{30}  &\blue{30}  &\blue{30}  &\blue{30}  &\blue{30}  &\blue{30}  \\
\hline
\end{tabular}
}\\
\vspace{1em}
\mysubtable{RMI, triangular mesh, $N_{\rm cells} = 16$}{
\centering
\begin{tabular}{|c||c|c|c|c|c|c|}
\hline
\diagbox{Solver}{Degree} & 1 & 2 & 3 & 4 & 5 & 6 \\
\hline
{Collocation} &\blue{30}  & \red{22.8}  & \red{21.52}  & \red{15.13}  & \red{8.841}  & \red{7.239}  \\
\hline
{Entropy projection} &\blue{30}  &\blue{30}  &\blue{30}  &\blue{30}  &\blue{30}  &\blue{30}   \\
\hline
\end{tabular}
}\\
\vspace{1em}
\mysubtable{RMI, triangular mesh, $N_{\rm cells} = 32$}{
\centering
\begin{tabular}{|c||c|c|c|c|c|c|}
\hline
\diagbox{Solver}{Degree} & 1 & 2 & 3 & 4 & 5 & 6 \\
\hline
{Collocation} & \blue{30}  & \red{23.84}  & \red{23.63}  & \red{8.752}  & \red{7.582}  & \red{3.946}  \\
\hline
{Entropy projection} &\blue{30}  &\blue{30}  &\blue{30}  &\blue{30}  &\blue{30}  &\blue{30}   \\
\hline
\end{tabular}
}\\
\vspace{.5em}
\caption{End time for simulations of the Richtmeyer-Meshkov instability on quadrilateral and triangular meshes.  On quadrilateral meshes, ``collocation'' refers to a nodal DGSEM discretization, while ``entropy projection'' refers to a method based on Gauss nodes.  On triangular meshes, ``collocation'' refers to nodal SBP discretization, while ``entropy projection'' refers to a modal entropy stable DG method.  Times colored \blue{blue} correspond to simulations which did not crash and ran to completion, while times colored \textcolor{red}{red} denote simulations which did crash.}
\label{tab:compare_rmi}
\end{table}

%\subsection{Three-dimensional experiments}

\subsubsection{Three-dimensional Kelvin-Helmholtz instability} For completeness, we also verify that a difference in robustness is observed for instability-type problems in three dimensions. Due to the high computational cost of entropy stable DG methods on tetrahedral meshes, we restrict ourselves to hexahedral meshes for the following experiments. We adapt the Kelvin-Helmholtz instability to three dimensions using the following initial condition:
\begin{align*}
&\rho = \frac{1}{2} + \frac{3}{4} B&
&p  = 1&\\
&u = \frac{1}{2} (B - 1)&
&v = \frac{1}{10} \sin(2 \pi x) \sin(2 \pi z)&
&w =\frac{1}{10} \sin(2 \pi x) \sin(2 \pi z),&
\end{align*}
where $B$ is defined as in (\ref{eq:B}). Table~\ref{tab:compare_khi3d} shows the results, which are similar to previous results for the two-dimensional test problems. We note for this example, both the relative and absolute adaptive time-step tolerances were set to $10^{-8}$ instead of $10^{-7}$. This was necessary to avoid crashes for the entropy projection method at degrees $N=6$ and $N=7$ on the finer $N_{\rm cells}=32$ mesh. 

\begin{table}
\centering
\mysubtable{3D KHI, $N_{\rm cells} = 16$}{
\centering
\begin{tabular}{|c||c|c|c|c|c|c|c|}
\hline
\diagbox{Solver}{Degree} & 1 & 2 & 3 & 4 & 5 & 6 & 7\\
\hline
{Collocation} & \blue{10}  & \red{2.73}  & \red{2.111}  & \red{1.978}  & \red{2.059}  & \red{1.797}  & \red{1.893}  \\
\hline
{Entropy projection} &\blue{10}  & \blue{10}  & \blue{10}  & \blue{10}  & \blue{10}  & \blue{10}  & \blue{10}  \\
\hline
\end{tabular}
}\\
\vspace{1em}
\mysubtable{3D KHI, $N_{\rm cells} = 32$}{
\centering
\begin{tabular}{|c||c|c|c|c|c|c|c|}
\hline
\diagbox{Solver}{Degree} & 1 & 2 & 3 & 4 & 5 & 6 & 7\\
\hline
{Collocation} & \red{4.049}  & \red{2.451}  & \red{2.061}  & \red{1.721}  & \red{2.071}  & \red{1.973}  & \red{1.952}  \\
\hline
{Entropy projection} &\blue{10}  & \blue{10}  & \blue{10}  & \blue{10}  & \blue{10}  & \blue{10}  & \blue{10}  \\
\hline
\end{tabular}
}\\
\caption{End time for simulations of the 3D Kelvin-Helmholtz instability on hexahedral meshes.  ``Collocation'' refers to a nodal DGSEM discretization, while ``entropy projection'' refers to a method based on Gauss nodes.  Times colored \blue{blue} correspond to simulations which did not crash and ran to completion, while times colored \textcolor{red}{red} denote simulations which did crash.}
\label{tab:compare_khi3d}
\end{table}

\subsection{Ideal GLM-MHD Equations}

Next, we consider the ideal GLM-MHD equations.
These equations use generalized Lagrange multiplier (GLM) technique to evolve towards a solution that bounds the magnetic field divergence. 
When the magnetic field divergence is non-zero, the GLM-MHD system requires the use of non-conservative terms to achieve entropy stability and to ensure Galilean invariance in the divergence cleaning technique.

The non-conservative GLM-MHD system without source terms reads
\begin{equation} \label{eq:noncons_sys}
\pd{\bm{u}}{t} + \nabla \cdot {\bm{f}(\bm{u})} + \bm{\Upsilon} = \bm{0},
\end{equation}
where the state variables are density, momentum, total energy, magnetic field, and the the so-called \textit{divergence-correcting field}, $\bm{u} = (\rho, \rho \fnt{v}, E, \fnt{B}, \psi)$, and the vectors $\fnt{v} = (u,v,w)$ and $\fnt{B} = (B_1,B_2,B_3)$ contain the velocities and magnetic field in $x$, $y$ and $z$, respectively.
The flux reads
\[
\bm{f}(\bm{u}) = 
\begin{pmatrix}
\rho \fnt{v} 
\\
\rho \fnt{v} \fnt{v}^T + \fnt{I} \left(p + \frac{1}{2} \|\fnt{B}\|^2 \right) - \fnt{B} \fnt{B}^T 
\\
\fnt{v} \left( \frac{1}{2} \rho \|\fnt{v}\|^2 + \frac{\gamma p}{\gamma-1} + \|\fnt{B}\|^2 \right) + \fnt{B} \left( {c_h} \psi - \left(\fnt{v}\cdot\fnt{B}\right) \right)
\\
\fnt{v} \fnt{B}^T - \fnt{B} \fnt{v}^{T} + \fnt{I} 
c_h \psi 
\\
c_h \fnt{B}
\end{pmatrix},
\]
where $\fnt{I}$ is again the $3 \times 3$ identity matrix, $\gamma$ is the heat capacity ratio, $c_h$ is the hyperbolic divergence-cleaning speed, and $p = (\gamma-1)\left(E - \left( \rho \|\fnt{v}\|^2 -  \|\fnt{B}\|^2 - \psi^2 \right)/2 \right)$ is the gas pressure.
Finally, the non-conservative term reads
\begin{align}
\bm{\Upsilon} &= \left(\nabla \cdot \fnt{B}\right) 
\LRp{ 0, \fnt{B}, \fnt{v} \cdot \fnt{B}, \fnt{v}, 0 } 
+
\LRp{ 0, \fnt{0}, \psi \left( \fnt{v} \cdot \nabla \psi \right), \fnt{0}, \fnt{v} \cdot \nabla \psi}.
\end{align}
To construct a two-dimensional version of the GLM-MHD system, we replace $\fnt{I}$ by a rectangular $3 \times 2$ identity matrix and neglect the flux in $z$.
However, we keep the third component of the velocity and magnetic field because plasma systems admit three-dimensional electromagnetic interactions in two-dimensional problems.
For details about the GLM-MHD system, we refer the reader to \cite{derigs2018ideal}.

The non-conservative GLM-MHD system \eqref{eq:noncons_sys} can be discretized using the collocation \eqref{eq:collocationDisc} and modal \eqref{eq:modalDisc} formulations by replacing the volume term $\fnt{F}$ and the surface term $\fnt{f}^*$ \cite{ruedaramirez2022entropy}.
In the collocation formulation the new terms read
\begin{gather}
\fnt{F}_{ij} = \bm{f}_S\LRp{{\fnt{u}}_i,{\fnt{u}}_j}
+ \bm{\Phi}^{\diamond}\LRp{{\fnt{u}}_i,{\fnt{u}}_j}, 
\qquad 
\fnt{f}^* = 
\begin{bmatrix}
\fnt{f}^*\LRp{{\fnt{u}}_1^+,{\fnt{u}}_1} 
+ \fnt{\Phi}^{\diamond}\LRp{{\fnt{u}}_1^+,{\fnt{u}}_1} \\
\fnt{f}^*\LRp{{\fnt{u}}_{N+1},{\fnt{u}}_{N+1}^+}
+
\fnt{\Phi}^{\diamond}\LRp{{\fnt{u}}_{N+1},{\fnt{u}}_{N+1}^+}
\end{bmatrix}.
\end{gather}
and in the modal formulation they read
\begin{gather}
\fnt{F}_{ij} = \bm{f}_S\LRp{\tilde{\fnt{u}}_i,\tilde{\fnt{u}}_j}
+ \bm{\Phi}^{\diamond}\LRp{\tilde{\fnt{u}}_i,\tilde{\fnt{u}}_j}, 
\qquad 
\fnt{f}^* = 
\begin{bmatrix}
\fnt{f}^*\LRp{\tilde{\fnt{u}}_1^+,\tilde{\fnt{u}}_1} 
+ \fnt{\Phi}^{\diamond}\LRp{\tilde{\fnt{u}}_1^+,\tilde{\fnt{u}}_1} \\
\fnt{f}^*\LRp{\tilde{\fnt{u}}_{N+1},\tilde{\fnt{u}}_{N+1}^+}
+
\fnt{\Phi}^{\diamond}\LRp{\tilde{\fnt{u}}_{N+1},\tilde{\fnt{u}}_{N+1}^+}
\end{bmatrix}.
\end{gather}
In addition to the symmetric two-point flux $\bm{f}_S$, we use a non-symmetric two-point term $\fnt{\Phi}^{\diamond}$ to account for the non-conservative term in the equation. The following experiment uses the non-conservative term presented by \citet{ruedaramirez2022entropy} and the entropy conservative flux of \citet{hindenlang2019new}, which is a natural extension of the entropy conservative, kinetic energy preserving, and pressure equilibrium preserving Euler flux of Ranocha \cite{ranocha2020entropy, ranocha2021preventing} to the GLM-MHD system.

\subsubsection{Two dimensional magnetized Kelvin-Helmholtz instability} 

To test the robustness of entropy projection schemes for the GLM-MHD system, we propose a modification of the Euler two-dimensional Kelvin-Helmholtz instability of Section~\ref{eq:KHI_euler}. 
The domain is $[-1,1]^2$ with the initial condition:
\begin{align}
&\rho = \frac{1}{2} + \frac{3}{4}B,&
&p = 1,& 
&\psi = 0,& \nonumber\\
&u = \frac{1}{2}(B-1),&
&v = \frac{1}{10}\sin\LRp{2\pi x},&
&w = 0,& \nonumber\\
&B_1 = 0,&
&B_2 = 0.125,&
&B_3 = 0,& \label{eq:mhdkhi_ic}
\end{align}
where $B(x,y)$ is as defined in (\ref{eq:B}). 
Each solver is run until final time $T_{\rm final} = 15$. 

For this example, we set $c_h$ as the maximum wave speed in the domain for the initial condition \eqref{eq:mhdkhi_ic} and keep it constant throughout the simulation.
This standard way of selecting $c_h$ has been shown to control the divergence error efficiently without affecting the time-step size \cite{derigs2018ideal, mignone2010high}.
We observed that smaller values of $c_h$ affect the robustness of the schemes for this problem, and higher values of $c_h$ increase the stiffness of the problem which can also lead to a crash if the tolerance for the adaptive time-stepping method is set too loosely. 

\begin{figure}
\centering
\includegraphics[width=\textwidth]{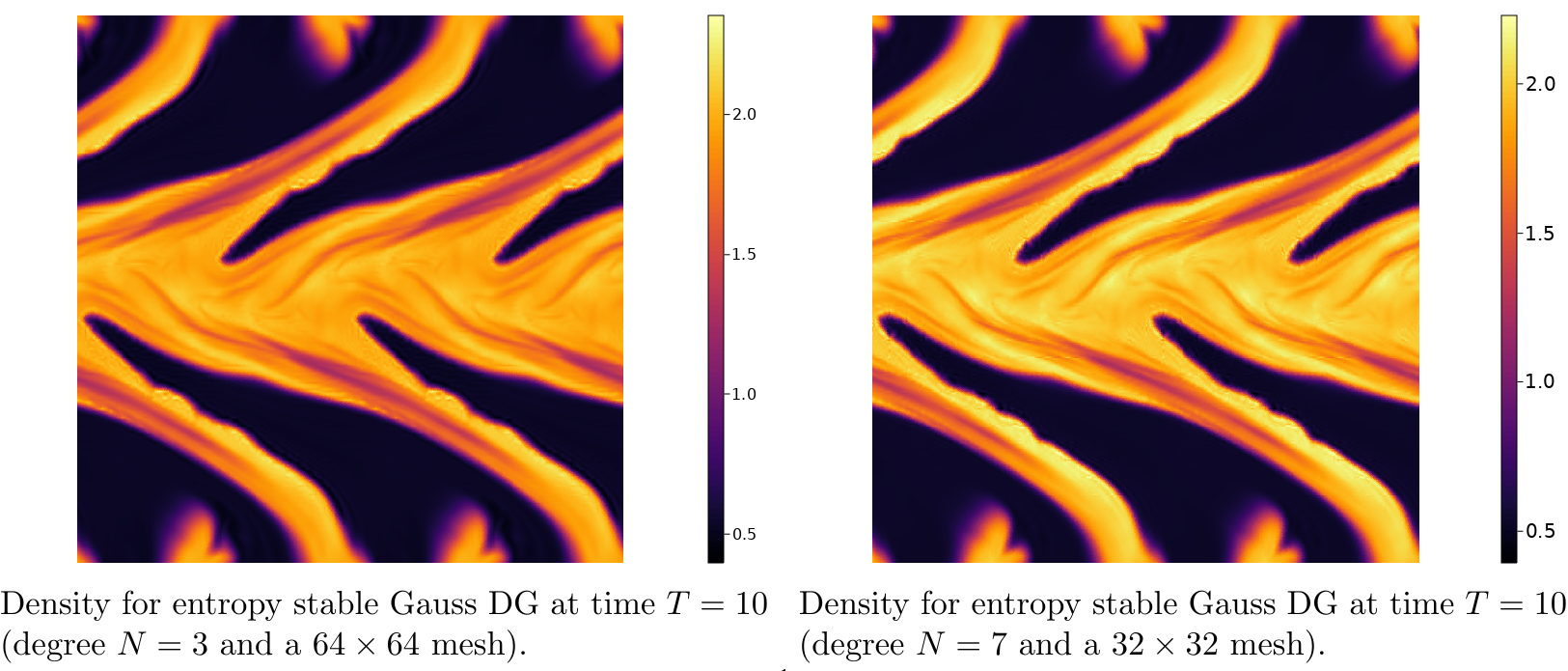}
%\subcaptionbox*{Density for entropy stable Gauss DG at time $T = 10$ (degree $N=3$ and a $64\times 64$ mesh).}{\includegraphics[width=.49\textwidth]{khi_mhd_p3_64_cells.png}}
%\hspace{.25em}
%\subcaptionbox*{Density for entropy stable Gauss DG at time $T = 10$ (degree $N=7$ and a $32\times 32$ mesh).}{\includegraphics[width=.49\textwidth]{khi_mhd_p7_32_cells.png}}
\caption{Snapshots of density for the magnetized Kelvin-Helmholtz instability using an entropy stable Gauss DG scheme on uniform quadrilateral meshes.}
\label{fig:khi_mhd_intro}
\end{figure}

Figure~\ref{fig:khi_mhd_intro} shows pseudocolor plots of the density at $T=10$ for the magnetized Kelvin-Helmholtz instability problem obtained with the entropy stable Gauss DG using polynomial degrees $N=3$ and $N=7$ on uniform meshes of $64 \times 64$ and $32 \times 32$ quadrilateral elements, respectively.
A comparison with Figure~\ref{fig:khi_intro} shows that the addition of a vertical magnetic field extends the flow features in the $y$ direction and suppresses many of the vortical structures at $T=10$. 
MHD turbulence eventually develops in the domain after $T=10$, which leads to later crash times for this example.

\begin{table}
\centering
\mysubtable{MHD KHI, quadrilateral mesh, $N_{\rm cells} = 16$}{
\centering
\begin{tabular}{|c||c|c|c|c|c|c|c|}
\hline
\diagbox{Solver}{Degree} & 1 & 2 & 3 & 4 & 5 & 6 & 7\\
\hline
{Collocation} & \blue{15} & \blue{15}  & \red{11.503}  & \red{10.988}  & \red{10.315}  & \red{10.230}  & \red{10.270}  \\
\hline
{Entropy projection} & \blue{15} & \blue{15} & \blue{15} & \blue{15} & \blue{15} & \blue{15} &  \blue{15} \\
\hline
\end{tabular}
}\\
\vspace{1em}
\mysubtable{MHD KHI, quadrilateral mesh, $N_{\rm cells} = 32$}{
\centering
\begin{tabular}{|c||c|c|c|c|c|c|c|}
\hline
\diagbox{Solver}{Degree} & 1 & 2 & 3 & 4 & 5 & 6 & 7\\
\hline
{Collocation} &  \blue{15} & \red{11.639} & \red{11.048} & \red{11.111} & \red{11.483} & \red{10.169} & \red{10.919} \\
\hline
{Entropy projection} & \blue{15} & \blue{15} & \blue{15} &  \blue{15} & \blue{15} & \blue{15} & \blue{15} \\
\hline
\end{tabular}
}\\
\vspace{1em}
\mysubtable{MHD KHI, triangular mesh, $N_{\rm cells} = 16$}{
\centering
\begin{tabular}{|c||c|c|c|c|c|c|}
\hline
\diagbox{Solver}{Degree} & 1 & 2 & 3 & 4 & 5 & 6 \\
\hline
{Collocation} & \red{12.846} & \red{13.797} & \red{10.626} & \red{10.212} & \red{10.990} & \red{9.973}\\
\hline
{Entropy projection} & \blue{15} & \blue{15} & \blue{15} & \blue{15} & \blue{15} & \blue{15} \\
\hline
\end{tabular}
}\\
\vspace{1em}
\mysubtable{MHD KHI, triangular mesh, $N_{\rm cells} = 32$}{
\centering
\begin{tabular}{|c||c|c|c|c|c|c|}
\hline
\diagbox{Solver}{Degree} & 1 & 2 & 3 & 4 & 5 & 6 \\
\hline
{Collocation} & \red{14.875} & \red{11.121} & \red{9.748} & \red{10.081} & \red{10.307} & \red{10.219}\\
\hline
{Entropy projection} &  \blue{15} & \blue{15}  & \blue{15}  & \blue{15}  & \blue{15}  & \blue{15} \\
\hline
\end{tabular}
}\\
\vspace{.5em}
\caption{End time for simulations of the magnetized Kelvin-Helmholtz instability on quadrilateral and triangular meshes. On quadrilateral meshes, ``collocation'' refers to a nodal DGSEM discretization, while ``entropy projection'' refers to a method based on Gauss nodes.  On triangular meshes, ``collocation'' refers to nodal SBP discretization, while ``entropy projection'' refers to a modal entropy stable DG method.  Times colored \blue{blue} correspond to simulations which did not crash and ran to completion, while times colored \textcolor{red}{red} denote simulations which did crash.}
\label{tab:compare_mhd_khi}
\end{table}

\subsection{Overview of results}

Tables~\ref{tab:compare_khi}, \ref{tab:compare_rti},  \ref{tab:compare_rmi} and \ref{tab:compare_mhd_khi} show what time the solver ran until for each solver on both quadrilateral and triangular meshes. We observe the pattern that, for degree $N > 1$, entropy stable methods which utilize the entropy projection appear be more robust than collocation-type schemes. Moreover, this pattern appears to hold independently of the polynomial degree and mesh size.

\subsection{Dependence of robustness on Atwood number}

While the numerical results in the previous section indicate a difference between different entropy stable schemes, they do not provide insight into why and when this difference in robustness manifests. The goal of this section is to establish a relationship between robustness, the Atwood number (a measure of the density contrast), and the use of the ``entropy projection'' in an entropy stable scheme. We restrict our focus to the Kelvin-Helmholtz instability for this section.

The results presented so far are somewhat unexpected, as the robustness of high order entropy stable DG schemes has been documented for a variety of flows where shocks and turbulent features are present \cite{gassner2016split, winters2018comparative, rojas2021robustness, parsani2021high}. In this section, we conjecture that the documented differences in robustness are due to the presence of both small-scale under-resolved features \textit{and} significant variations in the density. For example, entropy stable DGSEM methods are known to be very robust for the Taylor-Green vortex, where the density is near-constant throughout the duration of the simulation.

We examine the connection between density contrast and robustness by parametrizing the initial condition by the \textit{Atwood number}. Given a stratified fluid with two densities $\rho_1, \rho_2$, the Atwood number is defined as
\[
A = \frac{\rho_2-\rho_1}{\rho_1 + \rho_2} \in [0, 1),
\]
where it is assumed that $\rho_2 \geq \rho_1$. For a constant-density flow, $A = 0$, while $A\rightarrow 1$ indicates a flow with very large density contrasts. We investigate the behavior of different entropy stable methods for a version of the Kelvin-Helmholtz instability parametrized by the Atwood number $A$:
\begin{align*}
&\rho_1 = 1&
&\rho_2 = \rho_1 \frac{1+A}{1-A}&\\
%&\rho = \begin{cases}
%\rho_1 & y \not\in [-1/2, 1/2]\\
%\rho_2 & y \in [-1/2, 1/2]
%\end{cases}&
&\rho = \rho_1 + B (\rho_2 - \rho_1)&
&p = 1&\\
&u = B - \frac{1}{2}&
&v = \frac{1}{10}\sin(2\pi x)&
\end{align*}

\begin{figure}
\centering
\includegraphics[width=\textwidth]{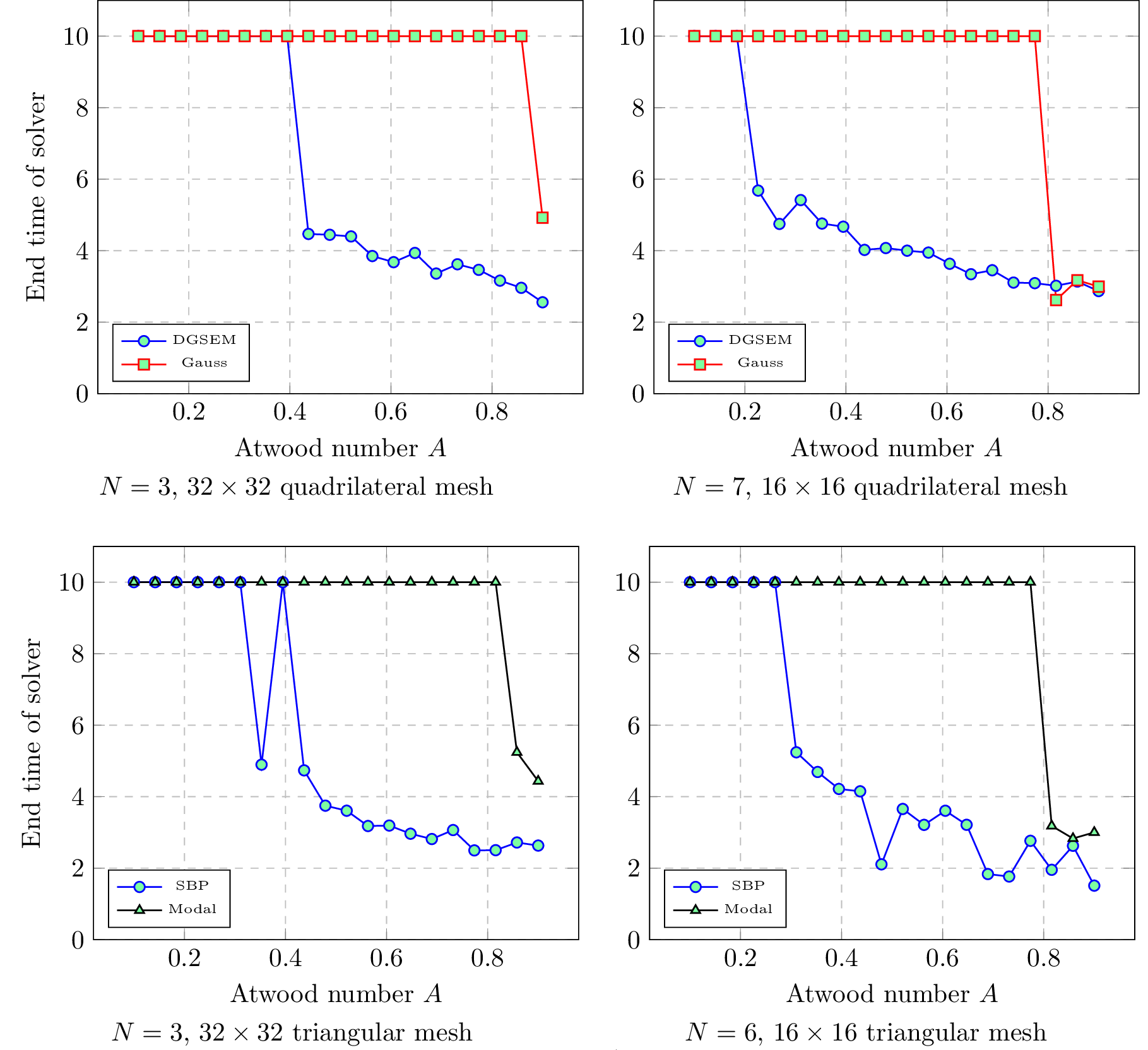}
\caption{Final times a solver an until as a function of Atwood number for the Kelvin-Helmholtz instability for DGSEM and various entropy stable solvers. End times less than final time $T_{\rm final} = 10$ indicate a crash.}
\label{fig:atwood}
\end{figure}

Figure~\ref{fig:atwood} shows the crash times for the Kelvin-Helmholtz instability using various entropy stable solvers at polynomial degrees $3$ and $7$. For quadrilateral meshes, we utilize entropy stable DGSEM solvers and entropy stable Gauss DG solvers. For triangular meshes, we utilize entropy stable multi-dimensional SBP solvers and entropy stable modal DG solvers. The DGSEM and SBP solvers are collocation-type schemes, while Gauss and modal DG solvers introduce the entropy projection.

For degree $3$ quadrilateral solvers, we utilize a $32\times 32$ mesh, while for degree $7$ quadrilateral solvers, we utilize a $16\times 16$ mesh. The mesh resolution is halved for polynomial degree $7$ simulations so that the total number of degrees of freedom is kept constant. For triangular solvers, we again use $32\times 32$ and $16\times 16$ uniform meshes, but we compare polynomial degrees $3$ and $6$, as SBP quadrature rules are available only up to degree $6$ in Trixi.jl. We run up to time $T_{\rm final}=10$ for $A \in [0.1, 0.9]$ and report the times each simulation ran until. For degree $N=3$, we observe that schemes which involve the entropy projection runs until the final time $T_{\rm final}=10$. Collocation-type schemes run to completion for low Atwood numbers, but crash earlier and earlier as the Atwood number increases. At degree $N=7$, we observe that while both collocation solvers and entropy projection solvers crash at higher Atwood numbers, entropy projection solvers begin crashing at higher Atwood numbers. For example, on quadrilateral meshes, DGSEM crashes around Atwood number $0.3$, while Gauss solvers crash around Atwood number $0.7$. We note that crash times for entropy projection schemes also tend to depend on the adaptive time-stepping tolerance. For example, for $N=3$ and a $32^2$ mesh, Gauss collocation runs stably to $T_{\rm final}=10$ if the absolute and relative tolerances are reduced to $10^{-9}$. The same is not true of entropy stable collocation-type schemes.

To provide another point of comparison, we ran simulations using an entropy stable DGSEM solver with sub-cell finite volume shock capturing \cite{hennemann2021provably} with Zhang-Shu positivity-preserving limiting for the density and pressure \cite{zhang2010positivity, zhang2012maximum}, which we refer to as DGSEM-SC-PP for shock capturing and positivity preservation.\footnote{For DGSEM-SC-PP, we utilize a 4-stage 3rd order adaptive strong stability preserving (SSP) Runge-Kutta time-stepping method \cite{kraaijevanger1991contractivity,conde2018embedded} with stepsize controller and efficient implementation of \cite{ranocha2021optimized}, which is necessary to ensure fully discrete positivity.} The entropy stable sub-cell finite volume-based shock capturing scheme utilizes a blending coefficient parameter $\alpha \leq \alpha_{\max}$ \cite{hennemann2021provably}. For these experiments, we set $\alpha_{\max} = 0.005$, which implies that the low order finite volume solution constitutes at most $0.5\%$ of the final blended solution. Despite the fact that this shock capturing is very weak, the resulting solver greatly improves robustness and enables long simulation times: for $N=3$ and a $32\times 32$ mesh, DGSEM-SC-PP runs stably to time $T_{\rm final} = 10$ for Atwood numbers up to $0.99$. However, we have also observed that the minimum value of $\alpha_{\max}$ necessary to avoid solver failure depends on the mesh resolution. For example, for $N=3$ and a $64\times 64$ mesh, we observe that DGSEM-SC-PP with $\alpha_{\max} = 0.005$ crashes around $t = 6.4871$.

\begin{remark}
We note that DGSEM with $\alpha_{\max} = 0.005$ shock capturing but no positivity preservation is not robust for the Kelvin-Helmholtz instability. For the initial condition (\ref{eq:khi_ic}), $N=3$, and a $64 \times 64$ mesh, DGSEM with shock capturing crashes around time $t=4.8891$. For $N=7$ and a $32\times 32$ mesh, DGSEM with shock capturing crashes around time $t=5.0569$. In contrast, DGSEM with only positivity preservation results in the simulation stalling due to a very small time-step. % DGSEM-SC crashes for Atwood numbers greater than $\approx 0.7$ for $N=3$ $\approx 0.6$ for $N=7$
\end{remark}

\section{The role of the entropy projection}
\label{sec:conjecture}

\subsection{Is robustness due only to the entropy projection?}

While the numerical results up to this point indicate that there is a significant difference in robustness for different entropy stable schemes, it is not yet clear that the increased robustness is due to the entropy projection. For example, the numerical experiments in Section~\ref{sec:comparison} compare entropy stable Gauss DG schemes to DGSEM on tensor product meshes and entropy stable ``modal'' DG methods to SBP schemes on triangular meshes. In both cases, a collocation scheme is compared to a scheme with higher accuracy numerical integration. Thus, it is not immediately clear whether the difference in robustness is due to the entropy projection or other factors such as the quadrature accuracy. We investigate whether the quadrature accuracy has a significant effect on stability by testing two additional variants of entropy stable DGSEM schemes on quadrilateral meshes. These schemes are purposefully constructed to be ``bad'' methods (in the sense that they introduce additional work without improving the expected accuracy), and are intended only to introduce the entropy projection. Both have quadrature accuracy similar to or lower than entropy stable DGSEM methods.

\begin{figure}
\centering
%\subcaptionbox*{Entropy stable Gauss DG}{\includegraphics[width=.32\textwidth]{khi_density_Gauss_p3_64_cells_t5.png}}%
%\hspace{.5em}%
%\subcaptionbox*{DGSEM with face-based entropy projection}{\includegraphics[width=.32\textwidth]{khi_density_hybrid_p3_64_cells_t5.png}}%
%\hspace{.5em}%
%\subcaptionbox*{DGSEM with volume-based entropy projection}{\includegraphics[width=.32\textwidth]{khi_density_staggered_p3_64_cells_t5.png}}%
\includegraphics[width=\textwidth]{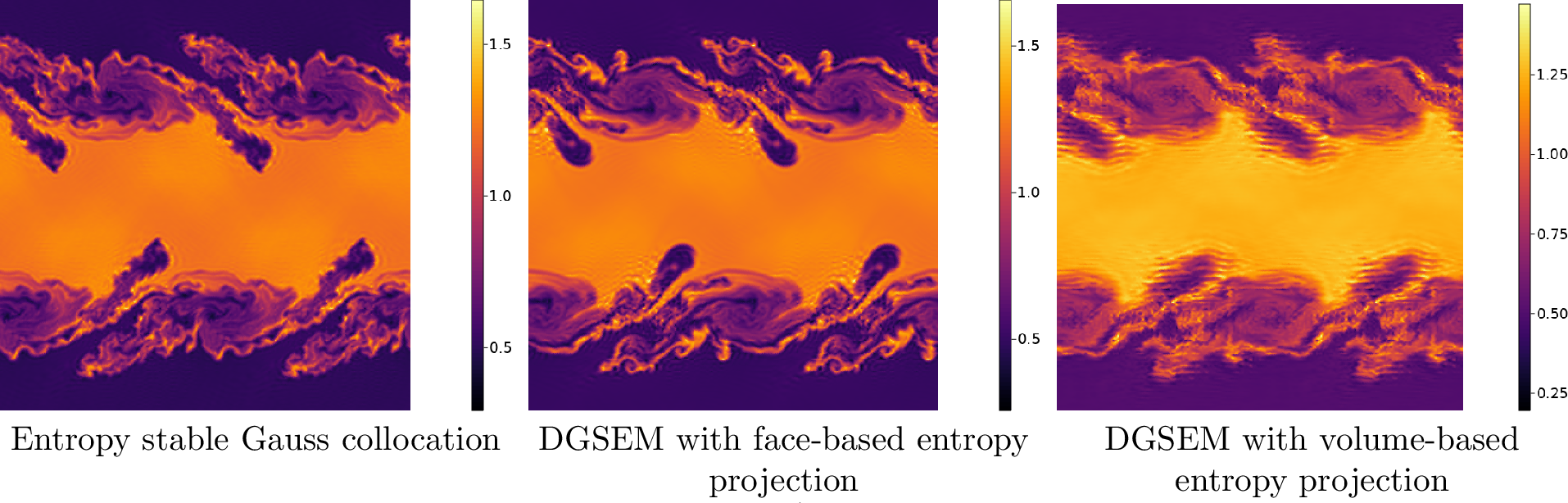}
\caption{Degree $N=3$ and $64\times 64$ grid Kelvin-Helmholtz simulations at $T=5$. All methods run until $T = 25$, while DGSEM crashes at $T \approx 3.5$.}%
\label{fig:ep_comparison}
\end{figure}

The first scheme utilizes LGL points for volume quadrature, but utilizes $(N+1)$ point Clenshaw-Curtis quadrature at the faces. This scheme can be directly derived from a modal formulation and (despite the lower polynomial exactness of Clenshaw-Curtis quadrature) can be shown to be entropy stable on affine quadrilateral meshes using the analysis in \cite{chan2019skew}. In order to retain entropy stability, the solution must be evaluated using the entropy projection at face nodes. We argue that the use of Clenshaw-Curtis quadrature does not result in a significant increase in quadrature accuracy over LGL quadrature: while Clenshaw-Curtis quadrature has been shown to be similar to Gauss quadrature for integration of analytic functions \cite{trefethen2008gauss}, for lower numbers of points we observe that the accuracy is comparable to LGL quadrature. Moreover, it was argued in \cite{kopriva2018stability} that increasing quadrature accuracy only for surface integrals or only for volume integrals does not provide sufficient anti-aliasing. We refer to this method as ``DGSEM with face-based entropy projection'' in Figure~\ref{fig:ep_comparison}.

\begin{remark}
We note that one can also build an entropy stable scheme from a combination of LGL volume points and Gauss face points. While this method possesses much of the simplicity and advantageous features of entropy stable DGSEM methods while also displaying improved robustness, this method results in a suboptimal rate of convergence by one degree \cite{chan2019skew}.
\end{remark}

The second scheme we test is similar to the staggered scheme of \cite{parsani2016entropy}. However, while the original scheme of Parsani et al.\ combines degree $N$ Gauss points with degree $(N+1)$ LGL points, we combine degree $N$ Gauss points with degree $N$ LGL points. This is a ``useless'' staggering in that it does not increase the accuracy of integration compared with DGSEM, and is intended only to introduce the entropy projection into the formulation.\footnote{This scheme can also be derived by beginning with an entropy stable DGSEM scheme and replacing the diagonal LGL mass matrix with the fully integrated dense mass matrix computed using Gauss quadrature. The resulting scheme can be made entropy stable by evaluating the spatial formulation using the entropy projection. More specifically, the appropriate entropy projection for this setting interpolates the entropy variables at Gauss nodes, then interpolates to LGL nodes.} We refer to this method as ``DGSEM with volume-based entropy projection'' in Figure~\ref{fig:ep_comparison}.

Figure~\ref{fig:ep_comparison} shows snapshots of the density for the Kelvin-Helmholtz instability for a degree $N=3$ mesh of $64 \times 64$ elements for each method. While the plots for the Gauss DG and DGSEM with face-based entropy projection have qualitative similarities, we observe that DGSEM with volume-based entropy projection results in a noisier solution. This may be due to inconsistency in terms of accuracy between the two quadrature rules used (e.g., $(N+1)$ point LGL and Gauss quadratures). However, all three entropy projection schemes remain stable, and we have verified that they are able to run until $T=25$ without crashing.

\begin{figure}
\centering
\includegraphics[width=\textwidth]{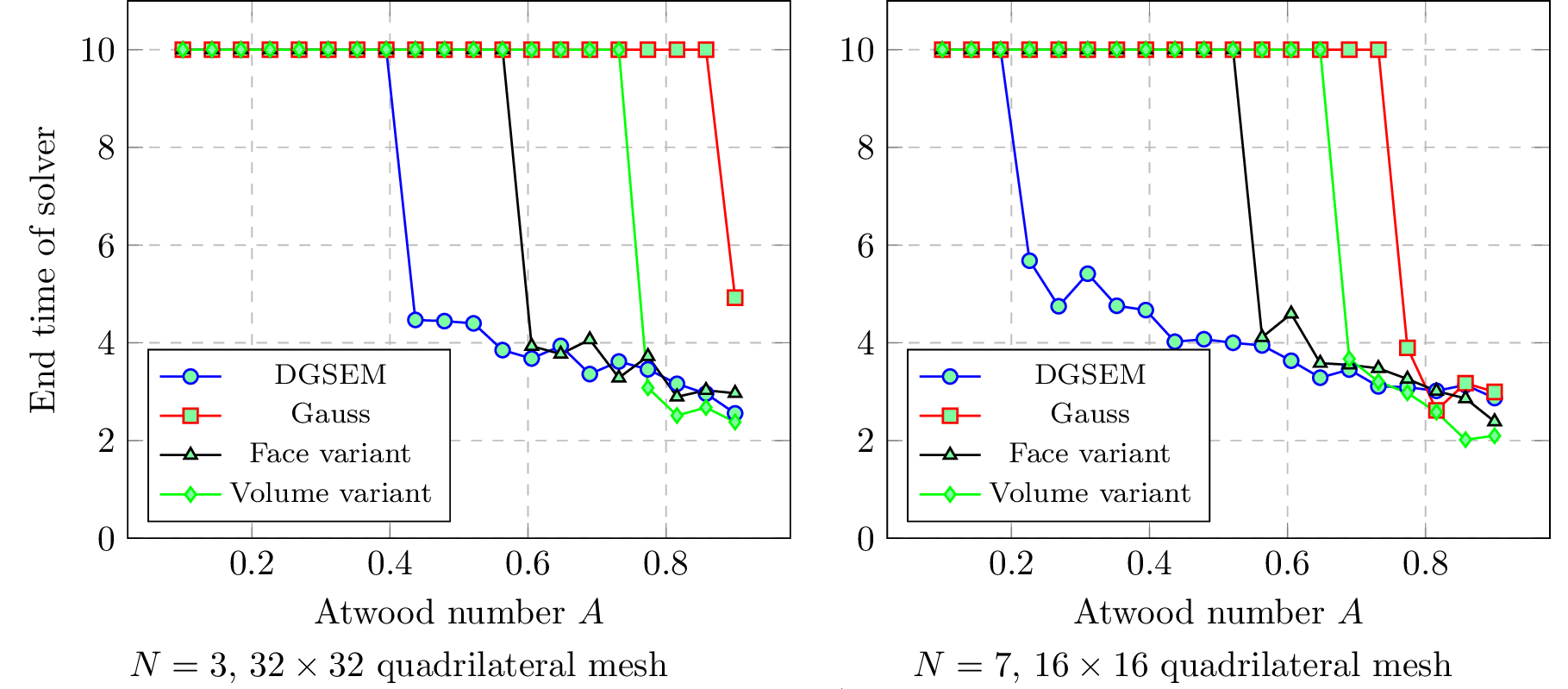}
\caption{Final times a solver an until as a function of Atwood number for the Kelvin-Helmholtz instability for DGSEM and different variants of entropy stable solvers based on the entropy projection. End times less than  final time $T_{\rm final} = 10$ indicate a crash.}
\label{fig:atwood2}
\end{figure}

We also compute crash times for each method for the Kelvin-Helmholtz instability with Atwood numbers $A \in [0.1, 0.9]$. These crash times are also compared to crash times of an entropy stable DGSEM method. These computations are performed on both a degree $N=3$ mesh of $32\times 32$ elements, as well as a degree $N=7$ mesh of $16 \times 16$ elements. Figure~\ref{fig:atwood2} plots the crash times for each method. We observe that all schemes which involve the entropy projection run stably for a wider range of Atwood numbers than entropy stable DGSEM, and that this effect becomes even more pronounced for degree $N=7$. However, for both the $N=3$ and $N=7$ experiments, the entropy stable Gauss schemes are stable for the widest ranges of Atwood numbers. 

These results indicate that incorporating the entropy projection does have a significant effect on the robustness of an entropy stable method, but that the entropy projection is not the only factor which impacts robustness. However, a detailed analysis of factors such as quadrature accuracy is out of the scope of this current work.

\subsection{Why is there a difference in robustness for different entropy stable methods?}

While the results from previous sections suggest that the entropy projection plays a role in the robustness of an entropy stable scheme, it is not clear \textit{why} it plays a role. While we do not have a thorough theoretical understanding of the entropy projection, initial experiments indicate that entropy projection schemes behave most differently from collocation schemes when the solution is either under-resolved or have near-zero density or pressure.

We illustrate the aforementioned behavior of the entropy projection using the one-dimensional compressible Euler equations. The conservative variables for the Euler equations are density, momentum, and total energy $(\rho, \rho u, E)$. Let $s(\bm{u}) = \log(p / \rho^\gamma)$ denote the specific entropy. The entropy variables for the convex entropy $S(\bm{u}) = -\rho s(\bm{u}) / (\gamma-1)$ are given by
\[
\bm{v}(\bm{u}) = \LRp{\frac{\gamma - s}{\gamma-1} - \frac{\rho u^2}{2p}, \frac{\rho u}{p}, -\frac{\rho}{p}}.
\]
Recall that the main steps of the entropy projection are as follows:
\begin{enumerate}
\item evaluate the entropy variables using degree $N$ polynomial approximations of the conservative variables
\item compute the quadrature-based $L^2$ projection of the entropy variables to degree $N$ polynomials
\item re-evaluate the conservative variables in terms of the projected entropy variables.
\end{enumerate}
These re-evaluated conservative variables are then used to compute contributions from an entropy stable DG formulation.

It was demonstrated numerically in \cite{chan2017discretely} that the entropy projection is high order accurate for sufficiently regular solutions. However, the behavior of the entropy projection was not explored for under-resolved or near-vacuum solution states. We illustrate this behavior using the following solution state:
\begin{gather}
\rho = 1 + e^{2 \sin(1 + k\pi x)}, \qquad u = \frac{1}{10} \cos(1 + k\pi x),\nonumber\\
p = p_{\min} + \frac{1}{2}\LRp{1 - \cos\LRp{k\pi x - \frac{1}{4}}},
\label{eq:entropy_projection_illustration}
\end{gather}
where $p_{\min} > 0$ is the minimum pressure, and $k$ is a parameter which controls the frequency of oscillation. As $k$ increases, the solution states in (\ref{eq:entropy_projection_illustration}) become more and more difficult to resolve, and as $p_{\min} \rightarrow 0$, the solution approaches vacuum and the entropy approaches non-convexity.
\begin{figure}
\centering
\includegraphics[width=.99\textwidth]{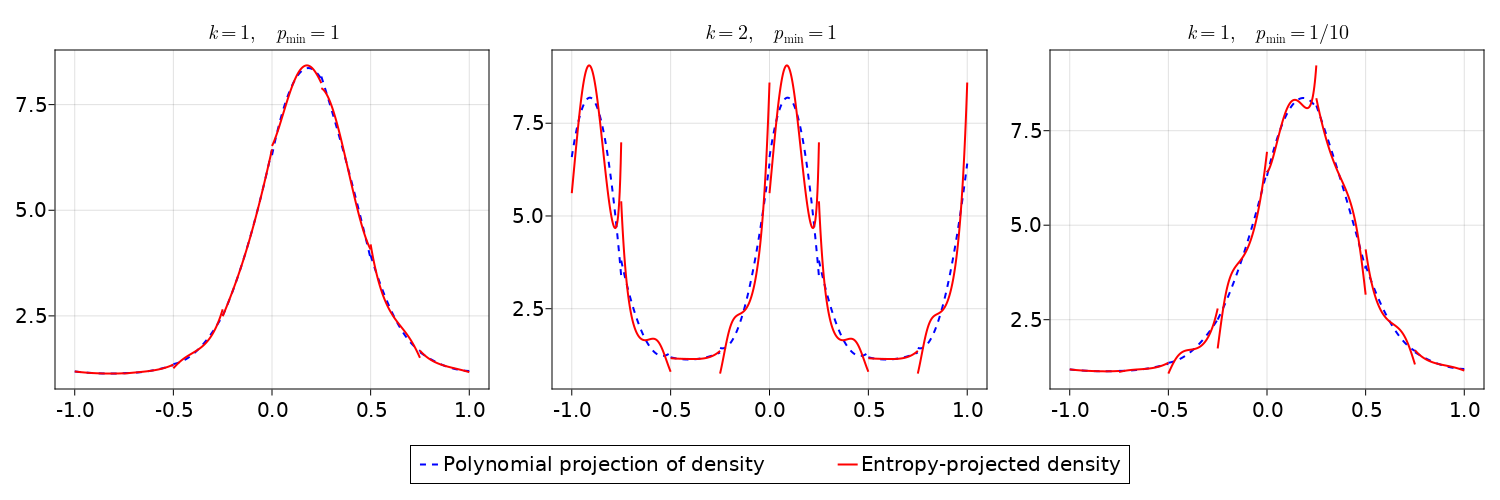}
\caption{Illustration of the effect of larger $k$ (under-resolution) and smaller $p_{\min}$ (near-vacuum state) on the entropy projection. A degree $N=2$ approximation and mesh of 8 elements were used. }
\label{fig:entropy_projection}
\end{figure}

Figure~\ref{fig:entropy_projection} illustrates the effect of increasing $k$ and decreasing $p_{\min}$ on the entropy projected conservative variables for a degree $N=2$ approximation on a coarse mesh of $8$ elements. As $k$ increases and the solution becomes under-resolved, the entropy projection develops large jumps at the interface. Similarly, as $p_{\min}$ decreases from $1$ to $1/10$, the entropy projection develops large jumps at the interface. We note that for both increased $k$ and decreased $p_{\min}$, spikes do not appear in the interior of the element.

%This phenomena may be related to the fact that the $L^2$ projection does not control \cite{houston2002discontinuous}

This indicates that the error in the entropy projection is influenced by both the numerical resolution and how close the entropy is to becoming non-convex. We denote the continuous entropy projection by $\tilde{\bm{u}} = \bm{u}\LRp{\Pi_N \bm{v}(\bm{u}_h)}$. Then, by the mean value theorem, we can bound the difference between the conservative and entropy-projected variables
\[
\nor{\bm{u}_h - \tilde{\bm{u}}}_{L^\infty} = \nor{\bm{u}_h - \bm{u}\LRp{\Pi_N\bm{v}(\bm{u}_h)}}_{L^\infty} \leq \nor{\pd{\bm{u}}{\bm{v}}}_{L^{\infty}} \nor{\bm{v}(\bm{u}_h) - \Pi_N\bm{v}(\bm{u}_h)}_{L^\infty},
\]
where $\pd{\bm{u}}{\bm{v}}$ is evaluated at some intermediate state between $\bm{u}_h$ and $\tilde{\bm{u}}$. The latter term in the bound $\nor{\bm{v}(\bm{u}_h) - \Pi_N\bm{v}(\bm{u}_h)}$ is small when the entropy variables are well-resolved, which we expect to be true when the solution is well-resolved and the mapping between conservative and entropy variables is well-conditioned. Conversely, high frequency components of the solution are often amplified when $\bm{v}(\bm{u})$ is highly nonlinear or the solution is under-resolved (this is the motivation behind filtering for stabilization \cite{orszag1971elimination, hesthaven2007nodal, bardos2015stability}), and we expect $\nor{\bm{v}(\bm{u}_h) - \Pi_N\bm{v}(\bm{u}_h)}$ to be large for such settings. The former term $\nor{\pd{\bm{u}}{\bm{v}}}$ is large when the mapping between conservative and entropy variables is nearly singular, which occurs when the entropy is nearly non-convex (for example, near-vacuum states).

\subsubsection{What role does entropy dissipation play?}

The previous section illustrates that entropy projection schemes are likely to differ from collocation schemes most when the solution is under-resolved or has near-zero density or pressure. Moreover, since the entropy projected variables in Figure~\ref{fig:entropy_projection} display spikes at the interfaces, it seems possible that the entropy projection would change the manner in which entropy dissipative interface dissipation terms are triggered. 
To test this hypothesis, we compute the evolution of entropy over time for the Kelvin-Helmholtz instability using both entropy stable Gauss DG and DGSEM-SC-PP, which is an entropy stable DGSEM with a shock capturing technique that consists in blending a sub-cell finite volume scheme with the DGSEM in an element-wise manner \cite{hennemann2021provably} and Zhang-Shu's positivity preserving limiter \cite{zhang2010positivity, zhang2012maximum}. The blending of the finite volume scheme is capped at $0.5\%$ in order to avoid unnecessary numerical dissipation.
We also compare entropy evolution for a scheme that blends a sub-cell finite volume scheme with the DGSEM in a subcell-wise manner \cite{rueda2022subcell}, which we refer to as DGSEM-subcell. 
The blending factors are chosen for each node (or subcell) to enforce lower bounds on density and pressure based on the low order solution, $\rho \ge 0.1 \, \rho^{\mathrm{FV}}$, $p \ge 0.1 \, p^{\mathrm{FV}}$. For this choice of lower bound, we observe high order accuracy for a two-dimensional sinusoidal entropy wave \cite{hindenlang2020order}. 
While this scheme is not provably entropy stable, it was demonstrated numerically in \cite{rueda2022subcell} that the use of subcell blending factors requires significantly lower levels of limiting compared with an element-wise limiting factor. 

%\note{TODO: Is DGSEM-subcell scheme high order accurate?} %TODO

\begin{figure}
\centering
\includegraphics[width=\textwidth]{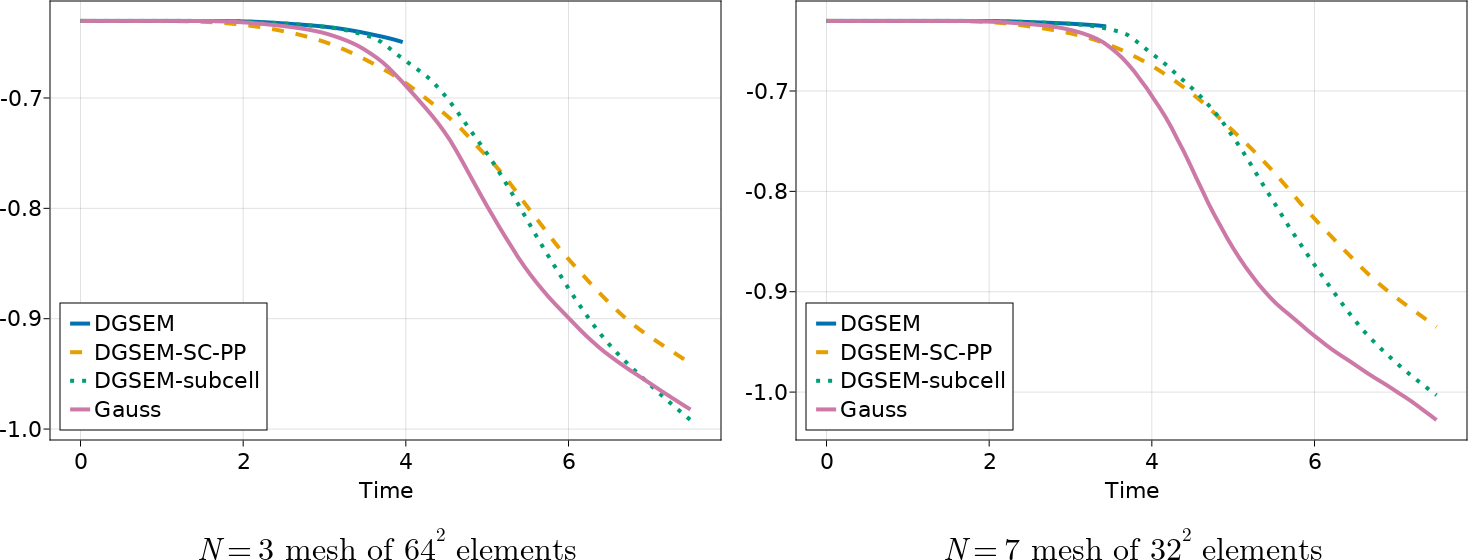}
%\subcaptionbox*{$N=3$ mesh of $64^2$ elements}{\includegraphics[width=.49\textwidth]{entropy_over_time_p3_64_cells.png}}
%\hspace{.1em}
%\subcaptionbox*{$N=7$ mesh of $32^2$ elements}{\includegraphics[width=.49\textwidth]{entropy_over_time_p7_32_cells.png}}
\caption{Evolution of entropy over time for the Kelvin-Helmholtz instability.}
\label{fig:entropy_evolution}
\end{figure}

Figure~\ref{fig:entropy_evolution} shows the evolution of the integrated entropy over the entire domain (which we have shifted to be positive) for the Kelvin-Helmholtz instability. Since periodic boundary conditions are used, the integrated entropy for the semi-discrete formulation can be proven to decrease over time. We observe that all four methods display similar entropy dissipation behavior until  time $t \approx 1.2$, after which DGSEM shows less entropy dissipation than either Gauss or DGSEM-SC-PP. However, while DGSEM-SC-PP initially dissipates more entropy than Gauss DG, the entropy dissipation for Gauss DG increases and overtakes that of DGSEM-SC-PP around time $t\approx 4$. Since entropy dissipation in both Gauss DG and DGSEM-SC-PP schemes is triggered by under-resolved flows (either through a modal indicator or through jump penalization terms) and since the Kelvin-Helmholtz instability generates increasingly small scales at larger times, this suggests that entropy dissipation for Gauss DG may be activated more strongly but at smaller scales than DGSEM-SC-PP. In contrast, Gauss DG dissipates more global entropy than DGSEM-subcell, though DGSEM-subcell eventually catches up to Gauss DG for $N=3$.

Our initial hypothesis was that the entropy projection in Gauss DG schemes results in larger interface jumps, which would trigger more entropy dissipation through jump penalization terms. However, this does not appear to be consistent with numerical results for entropy conservative schemes. To test these schemes, we focus on the three-dimensional Taylor-Green vortex. We note that the observed loss of robustness stands in stark contrast to the observed robustness of high order entropy stable and split-form DGSEM for the Taylor-Green vortex \cite{gassner2016split, chan2018efficient, rojas2021robustness}. This can be explained by the fact that the density remains near-constant over time for the Taylor-Green vortex; for a Kelvin-Helmholtz initial condition with a constant density, DGSEM runs stably up to final time $T=25$ for each of the previous numerical settings. Thus, while the Taylor-Green vortex generates small-scale flow features, it is a more benign test case when evaluating the robustness of high order entropy stable DG schemes.

However, when using a purely entropy conservative scheme (which can be constructed by utilizing entropy conservative interface fluxes), DGSEM methods can display non-robust behavior for the Taylor-Green vortex. We run the Taylor-Green vortex to final time $T_{\rm final} = 20$ using a variety of entropy conservative schemes: DGSEM, Gauss DG, as well as an entropy stable $C^0$ continuous Galerkin spectral element method (CGSEM) and a periodic finite difference method. We note that, because an entropy conservative scheme can be constructed given any summation-by-parts or skew-symmetric operator \cite{carpenter2014entropy, chen2017entropy, hicken2020entropy}, we are able to implement an entropy conservative $C^0$ continuous spectral element method and periodic finite difference method by constructing global difference operators from the tensor product of one-dimensional operators. These one-dimensional operators are provided by the Julia library SummationByPartsOperators.jl \cite{ranocha2021summationbypartsoperators}.

\begin{table}
\centering
\mysubtable{$N_{\rm cells} = 2^3$}{
\begin{tabular}{|c||c|c|c|c|c|c|c|}
\hline
\diagbox{Solver}{Degree} & 1 & 2 & 3 & 4 & 5 & 6 & 7\\
\hline
{DGSEM} & \blue{20}  & \blue{20}  & \blue{20}  & \blue{20}  & \red{16.4}  & \red{7.704}  & \red{7.482}  \\
\hline
{Gauss} & \blue{20}  & \blue{20}  & \blue{20}  & \blue{20}  & \blue{20}  & \blue{20}  & \blue{20}  \\
\hline
{CGSEM} & \blue{20}  & \blue{20}  & \blue{20}  & \blue{20}  & \blue{20}  & \blue{20}  & \blue{20}  \\
\hline
\end{tabular}
}
\\
\vspace{1em}
\mysubtable{$N_{\rm cells} = 4^3$}{
\begin{tabular}{|c||c|c|c|c|c|c|c|}
\hline
\diagbox{Solver}{Degree} & 1 & 2 & 3 & 4 & 5 & 6 & 7\\
\hline
{DGSEM} & \blue{20}  & \blue{20}  & \blue{20}  & \blue{20}  & \red{10.31}  & \red{5.792}  & \red{5.46}  \\
\hline
{Gauss} & \blue{20}  & \blue{20}  & \blue{20}  & \blue{20}  & \blue{20}  & \blue{20}  & \blue{20}  \\
\hline
{CGSEM} & \blue{20}  & \blue{20}  & \blue{20}  & \blue{20}  & \blue{20}  & \blue{20}  & \blue{20}  \\
\hline
\end{tabular}
}
\\
\vspace{1em}
\mysubtable{$N_{\rm cells} = 8^3$}{
\begin{tabular}{|c||c|c|c|c|c|c|c|}
\hline
\diagbox{Solver}{Degree} & 1 & 2 & 3 & 4 & 5 & 6 & 7\\
\hline
{DGSEM} & \blue{20}  & \blue{20}  & \blue{20}  & \blue{20}  & \red{6.035}  & \red{5.29}  & \red{5.02}  \\
\hline
{Gauss} & \blue{20}  & \blue{20}  & \blue{20}  & \blue{20}  & \blue{20}  & \blue{20}  & \blue{20}  \\
\hline
{CGSEM} & \blue{20}  & \blue{20}  & \blue{20}  & \blue{20}  & \blue{20}  & \blue{20}  & \red{17.5}  \\
\hline
\end{tabular}
}
\caption{End time for entropy \textit{conservative} simulations of the Taylor-Green vortex on hexahedral meshes. Times colored \blue{blue} correspond to simulations which did not crash and ran to completion, while times colored \textcolor{red}{red} denote simulations which did crash.}
\label{tab:tgv}
\end{table}

Table~\ref{tab:tgv} shows the end simulation time for each solver. We observe that again, despite the absence of any entropy dissipation, the Gauss DG solver is more robust than the DGSEM solver. The continuous spectral element solver CGSEM is also significantly more robust than the DGSEM solver, though it does lose robustness at higher orders and finer grid resolutions. We also ran periodic finite difference operators for grids with $4, 6, 8, 10, 12$ nodes in each dimension with orders of accuracy $2, 4, 6, 8, 10$. We observe that the periodic finite difference operator is as robust as the Gauss DG solver: for every grid resolution and order specified, the finite difference solver ran up to the final time $T_{\rm final} = 20$.

These experiments indicate that robustness for schemes involving the entropy projection is not solely due to the entropy dissipative terms. These experiments also show that robustness is improved for CGSEM and periodic finite difference solvers, neither of which contains interface terms. Since these results are on relatively coarse resolutions and utilize an entropy conservative scheme (when most practical schemes are entropy stable), further numerical experiments are necessary to carefully analyze the effect that different discretizations have on robustness.

\section{Applications toward under-resolved simulations}
\label{sec:spectra}

We conclude the paper with a discussion on a comparison between three schemes which include dissipative terms (entropy stable Gauss DG, entropy stable DGSEM-SC-PP, and DGSEM-subcell) for an under-resolved simulation. We run the Kelvin-Helmholtz instability using the initial condition (\ref{eq:khi_ic}), but modify the $y$-velocity perturbation to break symmetry of the resulting flow
\[
v = \frac{1}{10}\sin\LRp{2\pi x}\LRp{1 + \frac{1}{100}\sin\LRp{\pi x}\sin\LRp{\pi y}}.
\]
We run the simulation up to final time $T_{\rm final} = 25$. We use both a degree $N=3$ mesh of $64\times 64$ elements and a degree $N=7$ mesh of $32\times 32$ elements, each of which contains the same number of degrees of freedom. Due to the sensitivity of the Kelvin-Helmholtz instability problem and the long time window of the simulation, the results for each scheme are qualitatively very different.

\begin{figure}
\centering
\includegraphics[width=\linewidth]{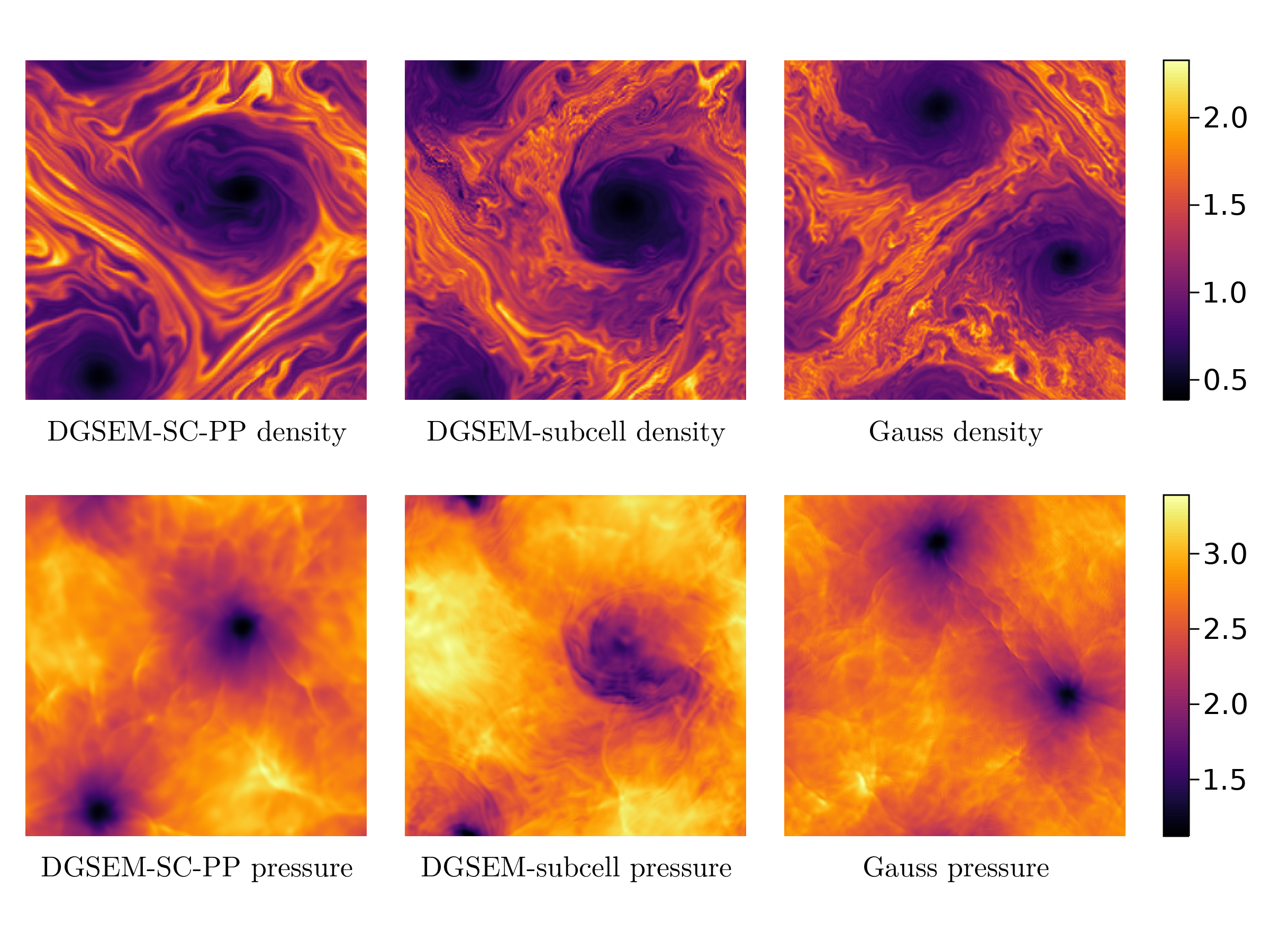}
%\subcaptionbox*{DGSEM-SC-PP density}{\includegraphics[width=.45\textwidth]{khi_density_DGSEM_SC_PP_p3_64_cells.png}}
%\hspace{1em}
%\subcaptionbox*{Gauss  density}{\includegraphics[width=.45\textwidth]{khi_density_Gauss_p3_64_cells.png}}\\
%\vspace{1em}
%\subcaptionbox*{DGSEM-SC-PP pressure}{\includegraphics[width=.45\textwidth]{khi_pressure_DGSEM_SC_PP_p3_64_cells.png}}
%\hspace{1em}
%\subcaptionbox*{Gauss pressure}{\includegraphics[width=.45\textwidth]{khi_pressure_Gauss_p3_64_cells.png}}
\caption{Density and pressure for the Kelvin-Helmholtz instability at $T_{\rm final}=25$ on a $N=3$ mesh of $64^2$ elements. }
\label{fig:khi_turbulence_N3}
\end{figure}

\begin{figure}
\centering
\includegraphics[width=\linewidth]{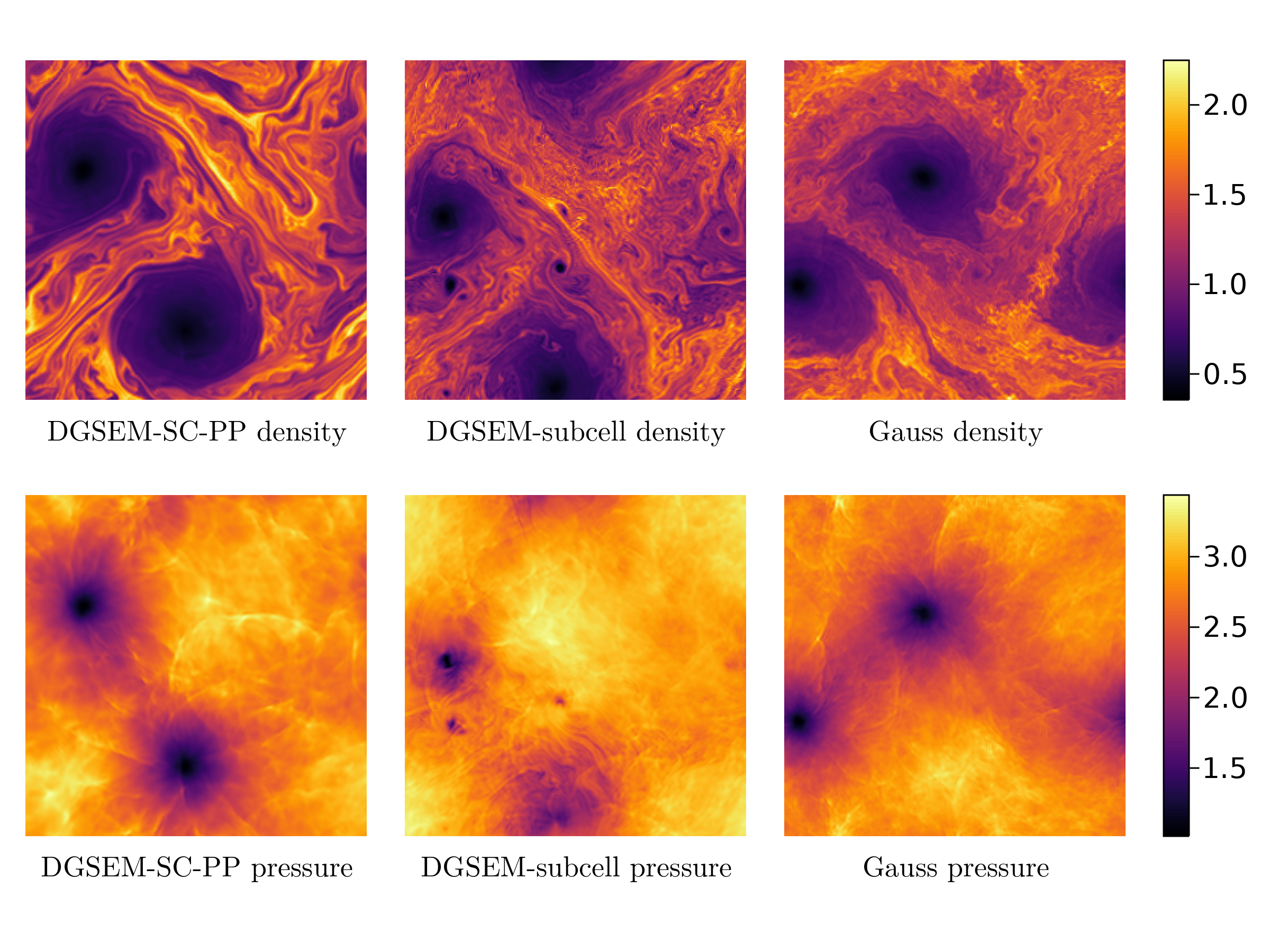}
%\subcaptionbox*{DGSEM-SC-PP density}{\includegraphics[width=.45\textwidth]{khi_density_DGSEM_SC_PP_p7_32_cells.png}}
%\hspace{1em}
%\subcaptionbox*{Gauss density}{\includegraphics[width=.45\textwidth]{khi_density_Gauss_p7_32_cells.png}}\\
%\vspace{1em}
%\subcaptionbox*{DGSEM-SC-PP pressure}{\includegraphics[width=.45\textwidth]{khi_pressure_DGSEM_SC_PP_p7_32_cells.png}}
%\hspace{1em}
%\subcaptionbox*{Gauss pressure}{\includegraphics[width=.45\textwidth]{khi_pressure_Gauss_p7_32_cells.png}}
\caption{Density and pressure for the Kelvin-Helmholtz instability at $T_{\rm final}=25$ on a $N=7$ mesh of $32^2$ elements. }
\label{fig:khi_turbulence_N7}
\end{figure}

Figures~\ref{fig:khi_turbulence_N3} and \ref{fig:khi_turbulence_N7} show snapshots of density and pressure for the entropy stable DGSEM-SC-PP and Gauss DG schemes. We observe that in both cases, the flow scales present in the DGSEM-SC-PP scheme are noticeably larger than those observed in the Gauss scheme. This is notable because the DGSEM-SC-PP scheme applies a very small amount of shock capturing: dissipation is added by blending the high order scheme with a low order finite volume scheme, and the amount of the blended low order solution is capped at $0.5\%$. However, even a small amount of dissipation produces a noticeable change on small-scale features in the resulting flow. We also observe the presence of shocklets or compression waves in the pressure, which mirror observations made in  \cite{terakado2014density}.\footnote{We note that these ``shocklets'' are not strictly shock waves, as the flow is not supersonic.}

\begin{figure}
\centering
\includegraphics[width=.95\textwidth]{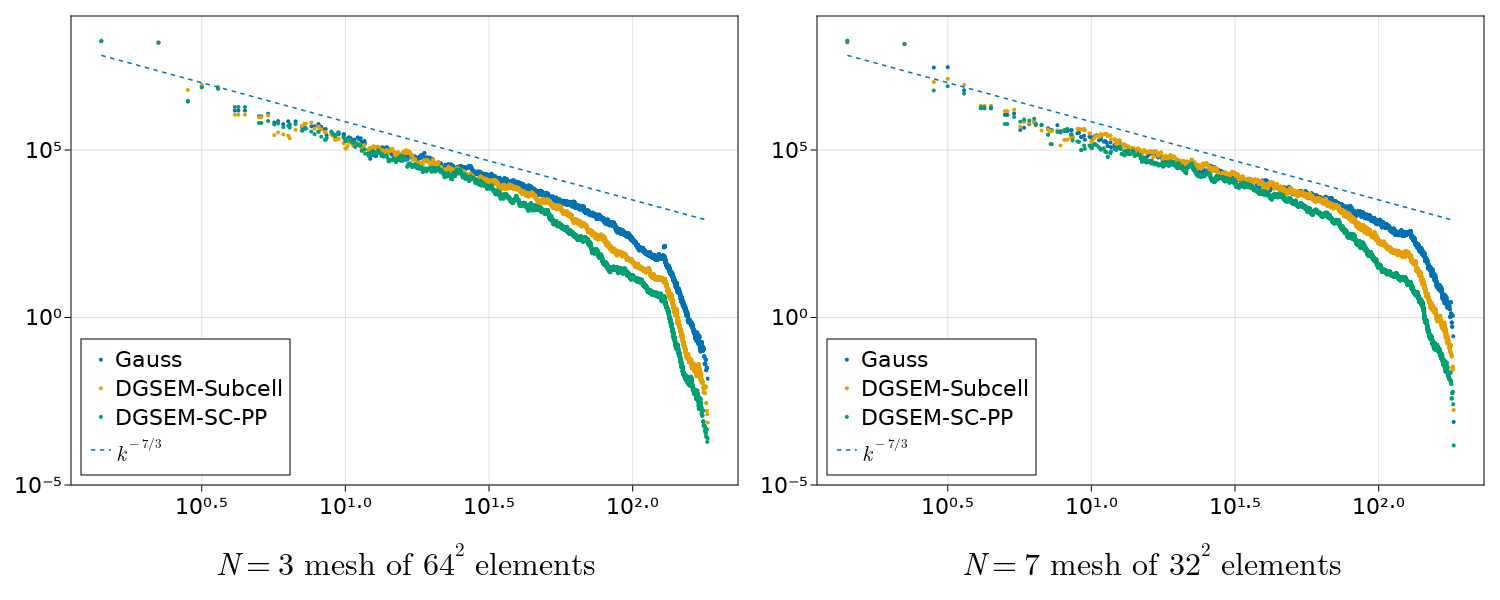}
%\subcaptionbox*{$N=3$, $64 \times 64$ mesh}{\includegraphics[width=.49\textwidth]{spectra_p3_64_cells_time_index_50.png}}
%\subcaptionbox*{$N=7$, $32 \times 32$ mesh}{\includegraphics[width=.49\textwidth]{spectra_p7_32_cells_time_index_50.png}}
\caption{Weighted power spectra for DGSEM-Subcell, entropy stable DGSEM-SC-PP, and entropy stable Gauss DG schemes. }%Left: $N=3$, $64\times 64$ mesh. Right: $N=7$, $32\times 32$ mesh.}
\label{fig:spectra}
\end{figure}

For $N=3$, the scales observed in DGSEM-subcell scheme are noticeably smaller than those of DGSEM-SC-PP but similar to those of the Gauss DG scheme. For $N=7$, the scales observed in the DGSEM-subcell scheme are again smaller than those of DGSEM-SC-PP, but appear to be slightly larger than those of the Gauss DG scheme. To avoid qualitative speculation, we compare these flows by computing the angle-averaged power spectra of the velocity weighted by $\sqrt{\rho}$ at final time $T_{\rm final}=25$ \cite{san2015evaluation, san2018stratified}. We follow \cite{moura2017eddy, winters2018comparative} and generate a grid of uniformly spaced points by evaluating the degree $N$ polynomial solution at $(N+1)$ equally spaced points along each dimension in the interior of each element of a uniform Cartesian mesh. The power spectra can then be computed from a fast Fourier transform of the resulting data. Figure~\ref{fig:spectra} shows the power spectra, which appear consistent with a $k^{-7/3}$ rate of decay from two-dimensional turbulence theory \cite{san2018stratified}. Moreover, we observe that the entropy stable Gauss DG scheme retains more energetic information than both DGSEM-SC-PP and DGSEM-subcell, though a spurious spike in the energy for Gauss DG schemes is observed near the higher wavenumbers for $N=3$.

\section{Conclusion}
\label{sec:conc}

This paper shows that for variable density flows which generate small-scale features, there are differences in robustness between entropy stable schemes which incorporate the entropy projection and those which do not. These differences in robustness are observed to depend on the Atwood number (measuring the density contrast) and persist across a range of polynomial degrees, mesh resolutions, and types of discretization. However, the mechanisms behind improved robustness for entropy projection schemes are currently unknown.
%\note{Finish conclusion.}

We note that any conclusions drawn concerning the robustness of DGSEM and Gauss DG should be restricted to the instability-type problems studied here. These results do not imply that Gauss is uniformly more robust than DGSEM.  %; one can construct problems with strong shocks where Gauss schemes behave less robustly than DGSEM.
Moreover, Gauss schemes are more computationally expensive than DGSEM schemes and result in smaller maximum stable timesteps \cite{gassner2011comparison, chan2015gpu, chan2018efficient, ranocha2021efficient}, so the appropriate scheme will depend on the use case. %utilizing a DGSEM scheme with nuanced shock capturing and positivity preservation may be more practical

\section*{Conflict of Interest Statement}

The authors declare that the research was conducted in the absence of any commercial or financial relationships that could be construed as a potential conflict of interest.

\section*{Author Contributions}

All authors contributed to the conception and design of the paper. JC, HR, and AR contributed numerical experiments. JC drafted the manuscript. All authors contributed to manuscript revision, read, and approved the submitted version. 
%The Author Contributions section is mandatory for all articles, including articles by sole authors. If an appropriate statement is not provided on submission, a standard one will be inserted during the production process. The Author Contributions statement must describe the contributions of individual authors referred to by their initials and, in doing so, all authors agree to be accountable for the content of the work. Please see \href{https://www.frontiersin.org/about/policies-and-publication-ethics#AuthorshipAuthorResponsibilities}{here} for full authorship criteria.

\section*{Funding}

Jesse Chan gratefully acknowledges support from the National Science Foundation under award DMS-CAREER-1943186.

Hendrik Ranocha was partially funded by the Deutsche Forschungsgemeinschaft
(DFG, German Research Foundation) under Germany's Excellence Strategy
EXC 2044-390685587, Mathematics M\"{u}nster: Dynamics-Geometry-Structure.

This work has received funding from the European Research Council through the
ERC Starting Grant ``An Exascale aware and Uncrashable Space-Time-Adaptive
Discontinuous Spectral Element Solver for Non-Linear Conservation Laws'' (Extreme),
ERC grant agreement no. 714487 (Gregor J.\ Gassner and Andr\'{e}s M. Rueda-Ram\'{i}rez).

Tim Warburton was supported in part by the Exascale Computing Project, a collaborative effort of two U.S. Department of Energy organizations (Office of Science and the National Nuclear Security Administration) responsible for the planning and preparation of a capable exascale ecosystem, including software, applications, hardware, advanced system engineering, and early testbed platforms, in support of the nation’s exascale computing imperative. Tim Warburton was also supported in part by the John K. Costain Faculty Chair in Science at Virginia Tech.

%Details of all funding sources should be provided, including grant numbers if applicable. Please ensure to add all necessary funding information, as after publication this is no longer possible.

\section*{Acknowledgments}
%This is a short text to acknowledge the contributions of specific colleagues, institutions, or agencies that aided the efforts of the authors.

%\note{Acknowledge Roci?}

This work used the Extreme Science and Engineering Discovery Environment (XSEDE) Expanse at the San Diego Supercomputer Center through allocation TG-MTH200014 \cite{towns2014xsede}.

This work was performed on the Cologne High Efficiency Operating Platform for
Sciences (CHEOPS) at the Regionales Rechenzentrum K\"{o}ln (RRZK) and on the
group cluster ODIN. We thank RRZK for the hosting and maintenance of the clusters.

The authors also thank Fabian F\"{o}ll for providing the initial condition for the Richtmeyer-Meshkov instability.

%\section*{Supplemental Data}
% \href{http://home.frontiersin.org/about/author-guidelines#SupplementaryMaterial}{Supplementary Material} should be uploaded separately on submission, if there are Supplementary Figures, please include the caption in the same file as the figure. LaTeX Supplementary Material templates can be found in the Frontiers LaTeX folder.

\section*{Data Availability Statement}
The code to generate and analyze the datasets used for this study can be found in the associated \href{https://github.com/trixi-framework/paper-2022-robustness-entropy-projection/}{reproducibility repository} \cite{reprorepo}. 
% Please see the availability of data guidelines for more information, at https://www.frontiersin.org/about/author-guidelines#AvailabilityofData

%\bibliographystyle{Frontiers-Harvard} %  Many Frontiers journals use the Harvard referencing system (Author-date), to find the style and resources for the journal you are submitting to: https://zendesk.frontiersin.org/hc/en-us/articles/360017860337-Frontiers-Reference-Styles-by-Journal. For Humanities and Social Sciences articles please include page numbers in the in-text citations
\bibliographystyle{Frontiers-Vancouver} % Many Frontiers journals use the numbered referencing system, to find the style and resources for the journal you are submitting to: https://zendesk.frontiersin.org/hc/en-us/articles/360017860337-Frontiers-Reference-Styles-by-Journal
\bibliography{refs}

\end{document}